\documentclass[11pt]{amsart}
\usepackage{float,url}
\usepackage{amssymb}
\theoremstyle{plain}
\newtheorem{thm}{Theorem}[section]
\newtheorem{lemma}[thm]{Lemma}
\newtheorem{cor}[thm]{Corollary}

\newtheorem{prop}[thm]{Proposition}

\theoremstyle{definition}

\newtheorem{rmk}[thm]{Remark}

\newcommand{\fra}{\mathfrak{a}}

\newcommand{\frg}{\mathfrak{g}}
\newcommand{\mfg}{\mathfrak{g}}
\newcommand{\frh}{\mathfrak{h}}

\newcommand{\frk}{\mathfrak{k}}
\newcommand{\mfk}{\mathfrak{k}}

\newcommand{\frp}{\mathfrak{p}}
\newcommand{\mfp}{\mathfrak{p}}

\newcommand{\frt}{\mathfrak{t}}
\newcommand{\mft}{\mathfrak{t}}

\newcommand{\bbR}{\mathbb{R}}

\newcommand{\bbZ}{\mathbb{Z}}

\newcommand{\C}{\mathcal{C}}

\newcommand{\F}{\mathcal{F}}

\newcommand{\caH}{\mathcal{H}}

\newcommand{\caM}{\mathcal{M}}
\newcommand{\M}{\mathcal{M}}
\newcommand{\caN}{\mathcal{N}}
\newcommand{\N}{\mathcal{N}}

\newcommand{\caS}{\mathcal{S}}

\usepackage{tikz-cd}

\begin{document}

\title[Enhanced Helgason-Johnson bounds for classical groups of type $B$]
{Enhanced Helgason-Johnson bounds for classical groups of type $B$}

\author{Zhan Ying}
\address[Ying]{School of Mathematical Sciences, Soochow University, Suzhou 215006,
P.~R.~China}
\email{zhanying\_1899@163.com}

\author{Chao-Ping Dong}
\address[Dong]{School of Mathematical Sciences, Soochow University, Suzhou 215006,
P.~R.~China}
\email{chaopindong@163.com}

\abstract{This paper improves the Helgason-Johnson bound in 1969 for non-Hermitian real classical groups of type $B$.}

\endabstract



\maketitle


\section{Introduction}

Let $G$ be a connected simple real non-compact Lie group with finite center. Let $\theta$ be the Cartan involution of $G$ so that $K:=G^{\theta}$ is a maximal compact subgroup of $G$.
Let $T$ be a maximal torus of $K$.  We have Cartan decomposition
$$
\frg_0=\frk_0 + \frp_0
$$
on the Lie algebra level. Let $\frg$ (resp. $\frk, \frt$) be the complexifications of $\frg_0$ (resp. $\frk_0, \frt_0$). By specifying a Vogan diagram of $\frg_0$, we have actually pinned down choices of positive roots
$\Delta^+(\frg, \frt)$ and $\Delta^+(\frk, \frt)$. The latter will be fixed once for all.
Let $\rho$ and $\rho_c$ be the half sum of roots in $\Delta^+(\frg, \frt)$ and $\Delta^+(\frk, \frt)$, respectively. By a \emph{$K$-type}, we mean an irreducible representation of $K$. We will denote a $K$-type by $E_{\eta}$, where $\eta$ is its highest weight living in $\mathcal{C}$, the dominant Weyl chamber of $\Delta^+(\frk, \frt)$.

Let $\pi$ be an irreducible admissible $(\frg, K)$ module. Take one of its lambda lowest $K$-types $E_{\mu}$. Let $G(\lambda_a(\mu))$ be the isotropic group of $\lambda_a(\mu)$ under the $G$ action. Let $H_{\rm s}=T_{\rm s} A_{\rm s}$ be a maximally split $\theta$-stable Cartan subgroup of $G(\lambda_a(\mu))$.
Then the infinitesimal character of  $\pi$ has the form
\begin{equation}\label{inf-char}
	\Lambda=\big(\lambda_a(\mu),\nu\big)\in \frh_{\rm s}^*=\frt_{\rm s}^*+\fra_{\rm s}^*.
\end{equation}
In representation theory of Lie groups, we care the most about the case that $\pi$ is unitary. Then the Helgason-Johnson bound obtained in \cite{HJ} says that
$$
\|\nu\|\leq \|\rho\|.
$$
The bound is attainable at the trivial representation. However, when $G$ is simple and non-compact, except for unitary characters, any irreducible unitary $(\frg, K)$ module is necessarily infinite-dimensional. Thus it is reasonable to assume that $\pi$ is irreducible, unitary and infinite-dimensional from now on. We will seek to improve the Helgason-Johnson bound under this reasonable assumption. Such an enhancement will help people to understand the unitary dual of $G$.

Since the main objects that we shall deal with are the lambda norm and spin norm of the $K$-types, one will agree that the hardest case is that $G$ is simply-connected: other groups are covered by the simply connected one and have fewer $K$-types. Thus we may and we will specify $G$ through its real Lie algebra $\frg_0$.

Recently, we have handled the $Sp(p, q)$ case \cite{DY2}, which happens to be connected and simply-connected. In the current paper, let us consider
\begin{equation}\label{g0}
\frg_0= \mathfrak{so}(2p, 2q+1), \quad 2\le p\leq q \, \mbox{  or  }\, p>q\ge 2
\end{equation}
Note that the above combines two classes of simple Lie algebras into one, see Appendix C of \cite{Kn}.
Also, to exclude the Hermitian symmetric groups,  we have cut off small values of $p$, $q$.

In this paper, we specify $\frg_0$ (thus the corresponding Lie groups) as in \eqref{g0}.
Firstly, assume that $2\le p\le q$ and let $n:=p+q$. Fixing an initial Vogan diagram for $\frg_0$  as in Appendix C of \cite{Kn}, we have
	\[\Delta^+(\mfg,\mft)=\{e_i\pm e_j\mid 1\le i<j\le n\}\cup\{e_t\mid 1\le t\le n\}\]
and
	\[\Delta^+(\mfk,\mft)=\{e_i\pm e_j\mid 1\le i<j\le p\}\cup \{e_{p+i}\pm e_{p+j}\mid 1\le i<j\le q\}\cup\{e_{p+j}\mid 1\le j\le q\}\]
and
	\[\Delta^+(\mfp,\mft)=\{e_i\pm e_{p+j}\mid 1\le i\le p\mid 1\le j\le q\}\cup \{e_i\mid 1\le i\le p\}.\]
Of course
\[\Delta^+(\mfg,\mft)=\Delta^+(\mfk,\mft)\cup\Delta^+(\mfp,\mft).\]
We use $\rho_n$ to denote half the sum of the elements of $\Delta^+(\mfp,\mft)$. Specifically,
\begin{align*}
	\rho&=\left(n-\frac{1}{2},n-\frac{3}{2},\dots,\frac{1}{2}\right),\\
	\rho_c&=\left(p-1, p-2, \dots,  0\ \left|\  q-\frac{1}{2},q-\frac{3}{2},\dots,\frac{1}{2}\right)\right.,\\
	\rho_n&=\left.\left(q+\frac{1}{2},\dots,q+\frac{1}{2}\ \right|\  0,\dots,0\right).
\end{align*}

Let $\C_\mfg$ be the dominant Weyl chamber corresponding to $\Delta^+(\mfg,\mft)$.  Define
\[W(\mfg,\mft)^1=\{w\in W(\mfg,\mft)\mid w(\C_\mfg)\subset \C\}.\]
Consider a sequence of real numbers $z_1,\dots,z_n$. Take out $p$ of its elements and arrange them in decreasing order as $u_1,\dots,u_p$. Arrange the remaining elements in their original relative order as $v_1,\dots,v_q$. For $w\in W(\mfg,\mft)^1$, we have
\[w(z_1,\dots,z_n)=(u_1,\dots,u_p\mid v_1,\dots,v_q)\quad\mbox{or}\quad w(z_1,\dots,z_n)=(u_1,\dots,-u_p\mid v_1,\dots,v_q).\]
Let $W(\mfg,\mft)^1_+$ collect the elements of $W(\mfg,\mft)^1$ whose action is merely a permutation of coordinates. For any $w\in W(\mfg,\mft)^1_+$, there is a counterpart of it in $W(\mfg,\mft)^1$ which changes the sign of $u_p$. Let us denote it by $w^-$. Thus
\[W(\mfg,\mft)^1=W(\mfg,\mft)^1_+\cup W(\mfg,\mft)^1_-,\]
where
\[W(\mfg,\mft)^1_-=\{w^-\mid w\in W(\mfg,\mft)^1_+\}.\]

Put $\rho^{(w)}:=w\rho$. Denote by $\rho_n^{(w)}$ the half sum of the roots in $w\Delta^+(\frp,\frt)$, then
\[\rho_n^{(w)}=\rho^{(w)}-\rho_c.\]

Given any $\eta\in\mathcal{C}$, we can find $w\in W(\frg, \frt)^1$ such that $\eta+2\rho_c\in \mathcal{C}_{w}$. As stated in \cite{SV}, there exists a unique
vector within $\mathcal{C}_{w}$ which is closest to $\eta+2\rho_c-\rho^{(w)}$. Denote it by $P(\eta+2\rho_c-\rho^{(w)})$. It turns out that it is independent of the choice of $w$. Let us denote it by $\lambda_a(\eta)$. Note that
$\|\lambda_a(\eta)\|$ is the \emph{lambda norm} of $\eta$ defined by Vogan in \cite{V81}.

Put
\begin{equation}\label{B0}
	B_0:=\left\{
	\begin{array}{ll}
		2pq(q-2)+\dfrac{4p^3+11p}{6}, & \mbox{if $q\ge p\ge 2$},\\
		(2q+1)\left(\bigg(p-\dfrac{3}{2}\bigg)^2+\dfrac{q(q+1)}{3}\right), & \mbox{if $p>q\ge 2$}.
	\end{array}
	\right.
\end{equation}
The meaning of $B_0$ will be explained in \eqref{B0-explanation}. Now let us state the main result.

\begin{thm}\label{thm-main}
Let $G$ be a connected Lie group with finite center whose Lie algebra $\frg_0$ is as in \eqref{g0}.
Let $\pi$ be any  infinite-dimensional irreducible unitary $(\frg, K)$ module whose
 infinitesimal character is given by \eqref{inf-char}. Then
	\[\Vert\nu\Vert\le \sqrt{B_0}.\]
\end{thm}

We believe that, as in the case of $Sp(p, q)$, this upper bound should be quite close to be tight.

The paper is organised as follows. Section \ref{sec-pre} collects necessary preliminaries. In particular, we recall spin norm, define $D_w(\mu)$ in \eqref{D-w-mu}, and discuss the structure of $\rho_n^{(w)}$.  In Sections \ref{sec-lambda-norm} and \ref{sec-Dw-mu}, for an arbitrary $G$, we obtain several equivalent forms of lambda norm, as well as the convexity of $D_w(\mu)$. Section \ref{sec-proof} proves Theorem \ref{thm-main}.

\medskip

\noindent \textbf{Notation.} Throughout this paper, we set
$\textbf{H}=\{\pm \frac{1}{2}, \pm \frac{3}{2}, \dots \}$ and $\textbf{I}=\{0, \pm 1, \pm 2, \dots\}$.
Write an $\Delta^+(\mfk,\mft)$-dominant weight as
\begin{equation}\label{mu}
	\mu=(a_1, \dots, a_p \mid b_1, \dots, b_q),
\end{equation}
where each entry lives in \textbf{H} or \textbf{I},  $a_1\ge a_2\ge \cdots \ge a_{p-1} \ge |a_p|$,
$b_1\ge b_2\ge \cdots b_q\ge 0$.
Moreover, we assume that $a_p\ge 0$. (We will show that this assumption involves no loss of generality.) We also write $a=(a_1,\dots,a_p)$ and $b=(b_1,\dots,b_q)$. Then $\mu=(a\mid b)$. If $a$ consists of half-integers and $b$ consists of integers, we write $a\in \textbf{H}$, $b\in \textbf{I}$ and $\mu \in (\textbf{H}\mid \textbf{I})$. Similar notation will also be used for the remaining cases. Set $P_i=p-i$ for $1\le i\le q$ and $Q_j=q-j+1/2$ for $1\le j\le q$. We use $\delta_{i,j}$ to denote the Kronecker delta, that is
\[\delta_{i,j}=\left\{
\begin{array}{ll}
	1, & \mbox{if $i=j$},\\
	0, & \mbox{otherwise}.
\end{array}
\right. \]
As a convention, when the lower limit of a summation exceeds the upper limit, the sum is understood to be zero.

\section{Preliminaries}\label{sec-pre}
Recall that we have fixed
$$
\Delta^+(\frg, \frt)=\Delta^+(\frk, \frt) \cup \Delta^+(\frp, \frt)
$$
in the introduction. Then $e_t-e_{t+1}(1\le t\le n-1)$ and $2e_n$ are the simple roots of $\Delta^+(\frg,\frt)$, with the corresponding fundamental dominant weights being $\xi_j=\sum_{t=1}^j e_t$, $1\le j\le n$. Note that $\beta:=e_1+e_{p+1}$ is the highest weight of $\frp$.

Let $\pi$ be an \emph{infinite dimensional} irreducible representation of $G$. Assume that its infinitesimal character is given by \eqref{inf-char}. Let $\mu$ be (the highest weight of) a $K$-type occurring in $\pi$. By \cite{V80},  the $K$-types of $\pi$ is a union of \emph{Vogan pencils}:
\begin{equation}
	\mathrm{Pencil}(\mu):=\{\mu+m\beta\mid m\in \bbZ_{\ge 0}\}.
\end{equation}
From \cite{D13},  we have the inequality
\begin{equation}\label{nu}
	\Vert \nu\Vert^2\le \min_{m\in\bbZ_{\ge 0}}\Vert \mu+m\beta\Vert_{\mathrm{spin}}^2-\Vert\mu\Vert_{\mathrm{lambda}}^2,
\end{equation}
where $\Vert\cdot\Vert_{\mathrm{spin}}$ and $\Vert\cdot\Vert_{\mathrm{lambda}}$ stand for spin norm and lambda norm, respectively. The notion spin norm will be introduced in \S\ref{sec-spin-norm} soon. For any $w\in W(\mfg,\mft)^1$, put
\begin{equation}\label{D-w-mu}
D_w(\mu):=\Vert\{\mu-\rho_n^{(w)}\}+\rho_c \Vert^2-\Vert\mu\Vert_{\mathrm{lambda}}^2
\end{equation}
and
\begin{equation}\label{B-mu}
B(\mu):=\min_{m\in\bbZ_{\ge 0}}\Vert \mu+m\beta\Vert_{\mathrm{spin}}^2-\Vert\mu\Vert_{\mathrm{lambda}}^2=\min_{m\in\bbZ_{\ge 0}}\min_{w\in W(\mfg,\mft)^1}D_w(\mu+m\beta).
\end{equation}
By \eqref{nu}, in order to prove Theorem \ref{thm-main}, we only need to show that $B(\mu)$ (or, more strongly, $D_w(\mu)$) satisfies the corresponding inequality.

\subsection{Spin norm}\label{sec-spin-norm}
Let $\mu$ be the highest weight of a $K$-type. Then the spin norm of $\mu$ is defined to be
 \begin{equation}\label{spin}
 	\Vert\mu\Vert_{\mathrm{spin}}:=\min_{w\in W(\frg,\frt)^1}\Vert \{\mu-\rho_n^{(w)}\}+\rho_c\Vert,
 \end{equation}
 where $\{\mu-\rho_n^{(w)}\}$ is the unique dominant weight that is conjugate to $\mu-\rho_n^{(w)}$ under the action of $W(\frk,\frt)$. Note that based on \cite{HP}, the spin norm was introduced in \cite{D13}. It offers a compact way of expressing the Dirac inequality.

Let  $x=(x_1,\dots,x_r)$ and $y=(y_1,\dots,y_r)$ be real vectors. Let $x_{(1)}\ge \cdots \ge x_{(r)}$ and $y_{(1)}\ge \cdots \ge y_{(r)}$ be the decreasing rearrangement of $x_1,\dots,x_r$ and $y_{1},\dots,y_r$, respectively. We say $x$ \emph{weakly majorizes} $y$ and write $x \succeq_w y$ if
\[\sum_{i=1}^{k}x_{(i)}\ge \sum_{i=1}^k y_{(i)}\]
for all $1\le k\le r$. We say that $x$  \emph{weakly dominates} $y$ and write $x\gg y$ if $x_i\ge y_i$ for all $1\le i\le r$. It is easy to see that if $x\gg y$, then $x\succeq_w y$.

Let  $u=(u_1,\dots,u_p)$ and $v=(v_1,\dots,v_q)$ be real vectors. Let $|u_{(1)}|\ge \cdots \ge |u_{(p)}|$ and $|v_{(1)}|\ge \cdots \ge |v_{(q)}|$ be the decreasing rearrangement of $|u_1|,\dots,|u_p|$ and $|v_{1}|,\dots,|v_q|$, respectively. Define
\[\Phi_p(u)=\sum_{i=1}^p(|u_{(i)}|+P_i)^2,\quad \Phi_q(v)=\sum_{j=1}^q(|v_{(j)}|+Q_j)^2,\]
where $P_i=p-i$ for $1\le i\le p$ and $Q_j=q-j+1/2$ for $1\le j\le q$.
\begin{lemma}\label{Phi}
	Let $r$ be $p$ or $q$. Let $u=(u_1,\dots,u_r)$ and $v=(v_1,\dots,v_r)$ be nonnegative $n$-dimensional real vectors. If $u$ weakly majorizes $v$, then
	\[\Phi_r(u)\ge \Phi_r(v).\]
\end{lemma}
\begin{proof}
	Note that $\Phi_p$ and $\Phi_q$ are increasing convex functions on the nonnegative real vectors. Then lemma immediately follows from Proposition B.2 in Chapter 4 of \cite{MI}.
\end{proof}
For convenience, for any $n$-dimensional vectors $\eta=(u_1,\dots,u_p\mid v_1,\dots,v_q)$, set
\begin{equation}\label{Phi-eta}
\Phi(\eta):=\Phi_p(u_1,\dots,u_p)+\Phi_q(v_1,\dots,v_q).
\end{equation}
Let $\nu_1$ and $\nu_2$ be any $\Delta^+(\mfk,\mft)$-dominant weights. We say that $\nu_1$ \emph{weakly majorizes} $\nu_2$ \emph{with respect to} $\Delta^+(\mfk,\mft)$, if the first $p$ components and the last $q$ components of $\nu_1$ weakly majorizes the corresponding components of $\nu_2$, respectively. For short, we write $\nu_1\succeq_\mfk \nu_2$. Then Lemma \ref{Phi} shows that if $\nu_1\succeq_\mfk \nu_2$, we have
\begin{equation}\label{compare-spin}
	\Phi(\nu_1)=\Phi_p(\nu_1)+\Phi_q(\nu_1)\ge \Phi_p(\nu_2)+\Phi_q(\nu_2)=\Phi(\nu_2).
\end{equation}
\subsection{The structure of $\rho_n^{(w)}$}
Define
\[\Omega_{p,q}^+=\{\rho_n^{(w)}\mid w\in W(\frg,\frt)^1_+\}.\]
Let us write a typical element in $\Omega_{p,q}^+$ as
\begin{equation}\label{rho-n-w}
\rho_n^{(w)}=(K_w\mid R_w)=\left.\left(k_1+\frac{1}{2},\dots,k_p+\frac{1}{2}\ \right|\ r_1,\dots,r_q\right).
\end{equation}
By convention, $k_0=q$,\ $r_0=p$ and $k_{p+1}=r_{q+1}=0$. For $1\le i\le q+1$, let
\[K_w(i):=(\underbrace{q-i+3/2,\dots,q-i+3/2}_{r_{i-1}-r_i})\]
if $r_{i-1}-r_i\neq 0$; let $K_w(i)$ be the empty array otherwise. For $1\le j\le p+1$, let
\[R_w(j):=(\underbrace{p-j+1,\dots,p-j+1}_{k_{j-1}-k_j})\]
if $k_{j-1}-k_j\neq 0$; let $R_w(j)$ be the empty array otherwise. Then
\[\rho_n^{(w)}=\big(K_w(1),\dots,K_w(q+1)\mid R_w(1),\dots,R_w(p+1)\big).\]

On the other hand, let
\[\Omega_{p,q}^-=\{\rho_n^{(w)}\mid w\in W(\mfg,\mft)^1_-\}.\]
The map
\[\mathrm{sgn}:\Omega_{p,q}^+\to \Omega_{p,q}^-,\quad \mathrm{sgn}(u_1,\dots,u_p\mid v_1,\dots,v_q)= (u_1,\dots,u_{p-1},-u_p\mid v_1,\dots,v_q)\]
is bijective. Note that $\mathrm{sgn}\in W(\mfg,\mft)$. Then
	\begin{align*}
	\Vert \mathrm{sgn}(\mu)\Vert_{\mathrm{spin}}&=\min_{w\in W(\mfg,\mft)^1}\Vert \{\mathrm{sgn}(\mu)-\rho_n^{(w)}\}+\rho_c\Vert\\
	&=\min_{w\in W(\mfg,\mft)^1}\Vert \{\mathrm{sgn}(\mu)-\mathrm{sgn}(\rho_n^{(w)})\}+\rho_c\Vert\\
	&=\min_{w\in W(\mfg,\mft)^1}\Vert \{\mu-\rho_n^{(w)}\}+\rho_c\Vert\\
	&=\Vert\mu\Vert_{\mathrm{spin}}.
\end{align*}
We will show that $\Vert\mathrm{sgn}(\mu)\Vert_{\mathrm{lambda}}=\Vert\mu\Vert_{\mathrm{lambda}}$ in Section \ref{sec-lambda-norm}. Therefore we can and will assume that the $p$-th component of $\mu$ is nonnegative.

\begin{prop}
	Let $\mu$ be any $\Delta^+(\mfk,\mft)$-dominant weight. Then the spin norm of $\mu$ can always be attained in $\Omega_{p,q}^+$.
\end{prop}
\begin{proof}
	Assume that the spin norm of $\mu$ is attained at $\rho_n^{(w_0)}\in \Omega_{p,q}^-$. Then $\mathrm{sgn}(\rho_n^{(w_0)})\in \Omega_{p,q}^+$. Recall that we have agreed that the $p$-th component of $\mu$ is nonnegative and the $p$-th component of $\rho_c$ is $0$. Then
	\[|\{\mu-\rho_n^{(w_0)}\}|+\rho_c\gg |\{\mu-\mathrm{sgn}(\rho_n^{(w_0)})\}|+\rho_c\]
	and
	\[\Vert|\{\mu-\rho_n^{(w)}\}|+\rho_c\Vert=\Vert\{\mu-\rho_n^{(w)}\}+\rho_c\Vert\]
	for any $w$ in $W(\mfg,\mft)^1$, where $|\eta|$ denotes the componentwise absolute value of $\eta$.  Lemma \ref{Phi} shows that the spin norm of $\mu$ is attained at $\mathrm{sgn}(\rho_n^{(w_0)})$ and we are done.
\end{proof}

Due to the above proposition, in this paper we consider only the case $\rho_n^{(w)}$ in $\Omega_{p,q}^+$. Furthermore, recall that $\Omega_{p,q}^*$ in $Sp(p,q)$ of \cite{DY2}, we have bijection
\begin{equation}\label{Omega-Sppq}
	\Omega_{p,q}^*\to \Omega_{p,q}^+,\quad (k_1,\dots,k_p\mid r_1,\dots,r_q)\mapsto \left.\left(k_1+\frac{1}{2},\dots,k_p+\frac{1}{2}\ \right| \ r_1,\dots,r_q\right).
\end{equation}

From now on, we determine $\rho_n^{(w)}:=(K_w\mid R_w)\in \Omega_{p,q}^+$ by choosing $(K_w^*\mid R_w^*)\in\Omega^*_{p,q}$. In such a situation, we always have
\[K_w=K_w^*+\bigg(\underbrace{\frac{1}{2},\dots,\frac{1}{2}}_p\bigg),\quad R_w=R_w^*.\]
	\begin{lemma}\label{Omega-pq}
	Let $(x_1,\dots,x_p\mid y_1,\dots,y_q)$ and $(k_1,\dots,k_p\mid r_1,\dots,r_q)$ be in $\Omega_{p,q}^+$. Then
	\[\sum_{i=1}^p |x_i-k_i|=\sum_{j=1}^q |y_j-r_j|.\]
\end{lemma}	
\begin{proof}
	Note that $|x_i-k_i|=|(x_i-1/2)-(k_i-1/2)|$, then the bijection \eqref{Omega-Sppq} and Lemma 2.6 of \cite{DY2} implies the lemma.
\end{proof}

\section{Lambda norm}\label{sec-lambda-norm}

\subsection{Some general results}

Recall that the lambda norm of a $K$-type has been defined in the introduction.
 This subsection aims to give some equivalent forms of the lambda norm for all real simple Lie groups.

Take $\mu$ as in \eqref{mu}. Choose $w\in W(\mfg,\mft)^1$ such that $\mu+2\rho_c\in \C_w$. We denote $w^{-1}(\mu+2\rho_c)-\rho$ by $\lambda^w(\mu)$ and denote $P_{C_\mfg}(\lambda^w(\mu))$ by $\lambda_a^w(\mu)$. Then
\begin{equation}\label{lambda-a}
	\lambda_a(\mu)=w\lambda_a^w(\mu).
\end{equation}

 \begin{prop}\label{lambda-rho}
 	Let $G$ be a connected simple real Lie group with finite center. Let $\mu$ be the highest weight of a $K$-type. Let $\C_\rho$ denotes the convex hull of the $W(\mfg,\mft)$ orbit of $\rho$. Then
 	\[\lambda_a(\mu)=\mu+2\rho_c-P_{\C_\rho}(\mu+2\rho_c),\]
 	where $P_{\C_w}$ and $P_{\C_\rho}$ are the projection operators onto $\C_w$ and $C_\rho$, respectively.	
 \end{prop}
 \begin{proof}
 	Choose $w\in W(\mfg,\mft)^1$ such that $\mu+2\rho_c\in \C_w$. By Corollary 1.9(a) of \cite{SV}, we have
 	\[\xi_\mu:=\mu+2\rho_c-\lambda_a(\mu)\in \mathrm{conv}(W_{\lambda_a(\mu)}\cdot \rho^{(w)})\subset \C_\rho,\]
 	where $W_{\lambda_a(\mu)}$ denotes the stabilizer of $\lambda_a(\mu)$ in $W(\mfg,\mft)$, while $\mathrm{conv}(W_{\lambda_a(\mu)}\cdot\rho^{(w)})$ denotes the convex hull of $W_{\lambda_a(\mu)}\cdot \rho^{(w)}$.

 We \emph{claim} that $\xi_\mu=P_{\C_\rho}(\mu+2\rho_c)$. By Theorem 3.16 of \cite{BC}, it suffices to show that for any $\eta\in \C_\rho$,
 	\begin{equation}\label{lambda-1}
 		\langle \mu+2\rho_c-\xi_\mu,\eta-\xi_\mu\rangle=\langle \lambda_a(\mu),\eta-\xi_\mu\rangle\le 0.
 	\end{equation}
 	Write $\eta$ as $\sum_{t}c_t w_t\cdot\rho$, where $c_t\ge 0$, $\sum_t c_t=1$, and $w_t\in W(\mfg,\mft)$. Then
 	\begin{equation}\label{lambda-2}
 		\begin{split}
 			\langle \lambda_a(\mu), \eta\rangle= \sum_t c_t \langle \lambda_a(\mu), w_t\cdot \rho\rangle &=\sum_t c_t \langle \lambda_a(\mu), w_tw^{-1}\cdot \rho^{(w)}\rangle\\
 			&\le \sum_t c_t \langle \lambda_a(\mu), \rho^{(w)}\rangle=\langle \lambda_a(\mu), \rho^{(w)}\rangle.
 		\end{split}
 	\end{equation}
 	On the other hand, write $\xi_\mu$ as $\sum_{s}d_s w_s\cdot\rho^{(w)}$, where $d_s\ge 0$, $\sum_s d_s=1$, and $w_s\in W_{\lambda_a(\mu)}$. Then $w_s\lambda_a(\mu)=\lambda_a(\mu)$ shows that
 	\begin{equation}\label{lambda-3}
 		\langle \lambda_a(\mu), \xi_\mu\rangle=\sum_{s}d_s \langle \lambda_a(\mu), w_s\cdot\rho^{(w)}\rangle=\sum_{s}d_s \langle \lambda_a(\mu), \rho^{(w)}\rangle=\langle \lambda_a(\mu), \rho^{(w)}\rangle.
 	\end{equation}
 	Combining \eqref{lambda-2} and \eqref{lambda-3} gives \eqref{lambda-1}. Therefore the claim holds and we are done.
 \end{proof}

\subsection{Properties of $\mathfrak{so}(2p, 2q+1)$}

From now on, let us focus on those $G$ whose Lie algebra is given by \eqref{g0}.
Then
$$
\C_\mfg=\{(z_1,\dots,z_n) \mid z_1\ge\cdots\ge z_n\ge 0\}.
$$

\begin{lemma}\label{Weyl-B}
Let $u$ and $v$ be any $n$-dimensional vectors. We have
	\begin{itemize}
		\item[(a)] If $u\gg v$, then $P_{C_\mfg}(u)\gg P_{C_\mfg}(v)$.
		\item[(b)] $P_{C_\mfg}(u)$ weakly majorizes $u$.
	\end{itemize}
\end{lemma}
\begin{proof} For (a), let $P_{\C_\mfg}(u)=(x_1,\dots,x_n)$ and $P_{\C_\mfg}(v)=(y_1,\dots,y_n)$. Set
	\[\eta_{\mathrm{min}}=(\min\{x_1,y_1\},\dots,\min\{x_n,y_n\}),\quad \eta_{\mathrm{max}}=(\max\{x_1,y_1\},\dots,\max\{x_n,y_n\}).\]
	Then $\eta_{\mathrm{min}}$ and $\eta_{\mathrm{max}}$ are in $\C_\mfg$. By the definition of the projection operator,
	\[\Vert u-P_{\C_\mfg}(u)\Vert\le\Vert u-\eta_{\mathrm{max}}\Vert,\quad \Vert v-P_{\C_\mfg}(v)\Vert\le\Vert v-\eta_{\mathrm{min}}\Vert.\]
	Then
	\[\Vert u-P_{\C_\mfg}(u)\Vert^2-\Vert u-\eta_{\mathrm{max}}\Vert^2+\Vert v-P_{\C_\mfg}(v)\Vert^2-\Vert v-\eta_{\mathrm{min}}\Vert^2\le 0.\]
	On the other hand, we have
	\begin{align*}
		&\quad\Vert u-P_{\C_\mfg}(u)\Vert^2-\Vert u-\eta_{\mathrm{max}}\Vert^2+\Vert v-P_{\C_\mfg}(v)\Vert^2-\Vert v-\eta_{\mathrm{min}}\Vert^2\\
		&=2\langle u,\eta_{\mathrm{max}}-P_{C_\mfg}(u)\rangle+2\langle v,\eta_{\mathrm{min}}-P_{C_\mfg}(v)\rangle\\
		&\ge 2\langle u,\eta_{\mathrm{max}}-P_{C_\mfg}(u)\rangle+2\langle u,\eta_{\mathrm{min}}-P_{C_\mfg}(v)\rangle\\
		&=2\langle u,0\rangle=0,
	\end{align*}
where the first step uses $\Vert \eta_{\max}\Vert^2+\Vert \eta_{\mathrm{min}}\Vert^2=\Vert P_{C_\mfg}(u)\Vert^2+\Vert P_{C_\mfg}(v)\Vert^2$, and the second step uses $\langle u-v,P_{C_\mfg}(v)-\eta_{\mathrm{min}} \rangle\ge 0$. Thus
	\[\Vert u-P_{\C_\mfg}(u)\Vert=\Vert u-\eta_{\mathrm{max}}\Vert.\]
	Therefore, the uniqueness of the projection forces $P_{C_\mfg}(u)=\eta_{\mathrm{max}}$. This proves (a).
	
	For (b), note that $P_{\C_\mfg}(u)+\sum_{i=1}^re_i$ is in $\C_\mfg$ for any $1\le r\le n$. Then Theorem 3.16 of \cite{BC} shows that
	\[\langle u-P_{\C_\mfg}(u), \sum_{i=1}^re_i\rangle\le 0.\]
	This is equivalent to (b).
\end{proof}

\begin{cor}[Monotonicity of the lambda norm]\label{lambda-up}
Let $\nu$ and $\eta$ be any $\Delta^+(\mfk, \mft)$-dominant weights. If $\nu\gg \eta$, then $\Vert \nu\Vert_{\mathrm{lambda}}\ge \Vert \eta\Vert_{\mathrm{lambda}}$.
\end{cor}
\begin{proof}
 Note that $\nu\gg\eta$  impiles that $\lambda^{w_1}(\nu)\gg \lambda^{w_2}(\eta)$.  Thus Lemma \ref{Weyl-B}(a) shows that $\lambda_a^{w_1}(\nu)\gg \lambda_a^{w_2}(\eta)$. Therefore,
	\[\Vert \nu\Vert_{\mathrm{lambda}}=\Vert \lambda_a^{w_1}(\nu)\Vert\ge \Vert \lambda_a^{w_2}(\eta)\Vert=\Vert \eta\Vert_{\mathrm{lambda}}.\]
\end{proof}

Let us give an explicit formula for lambda norm. Before that, we need some notation. Write  $w^{-1}(\mu+2\rho_c)$ as $(u_1,\dots,u_n)$. Define
\begin{align*}
	\M(\mu)&=\mathrm{card}\{1\le i\le n-1\mid u_i=u_{i+1}\},\\
	\mathcal{N}(\mu)&=\mathrm{card}\{1\le i\le n\mid \mbox{$u_i=n-i$ and $u_{i}\neq u_{i+1}$}\}
\end{align*}
and
\[\mathcal{H}(\mu)=\mathrm{card}\{1\le i\le n-1\mid u_i-u_{i+1}=\frac{1}{2}\ \mbox{and}\ u_i\neq n-i\},\]
where we adopt the convention that $u_n\neq u_{n+1}$. When the components of $\mu$ are either all integers or all half-integers, we have
\[\N(\mu)=\left\{
\begin{array}{ll}
	2c_2+1,& \mbox{if}\ c_1>c_2,\\
	2c_1, & \mbox{if}\ c_1\le c_2,
\end{array}
\right.\]
where $c_1$ and $c_2$ denote the least indices of zero components in $a$ and $b$, respectively.
\begin{prop}\label{lambda norm}
	 Take $\mu$ as in \eqref{mu}. Then
	\begin{equation}\label{lambda norm-strat}
		\Vert\mu\Vert_{\mathrm{lambda}}^2=\Vert \lambda^w(\mu)\Vert^2-\frac{\M(\mu)}{2}-\frac{\N(\mu)}{4}-\frac{\caH(\mu)}{8}.
	\end{equation}
\end{prop}

\begin{proof}
 Recall
	\[C_\mfg=\left\{(z_1,\dots,z_n)\in\bbR^n\mid z_1\ge\cdots\ge z_n\ge 0\right\}.\]
Note that in both the first $p$ and last $q$ components of
$\mu+2\rho_c$, the difference between any two adjacent components is at least $2$. We will compute $\lambda_a^w(\mu)$, thereby obtaining lambda norm from \eqref{lambda-a}.

	Write $w^{-1}(\mu+2\rho_c)$ and $\lambda^w(\mu)$ as $(u_1,\dots,u_n)$ and $(t_1,\dots,t_n)$, respectively. Put $N_i=n-i+\frac{1}{2}$. For $1\le i\le n-1$,
	\[t_i-t_{i+1}=(u_i-N_i)-(u_{i+1}-N_{i+1})=(u_i-u_{i+1})-(N_i-N_{i+1})=(u_i-u_{i+1})-1.\]
	So $t_i< t_{i+1}$ if and only if $u_i=u_{i+1}$ or $u_i=u_{i+1}+1/2$. On the other hand, let $\rho_c^*$ be the rearrangement of $\rho_c$ in decreasing order. Then $w^{-1}(\mu+2\rho_c)$ weakly dominates $2\rho_c^*$. Note that
$$
2\rho_c^*+\left(\frac{1}{2},\dots,\frac{1}{2}\ \left|\ \frac{1}{2},\dots,\frac{1}{2}\right)\right.
$$ weakly dominates $\rho$. This shows that the only possible negative component of $\lambda^w(\mu)=w^{-1}(\mu+2\rho_c)-\rho$ is $-1/2$.
	
We have the following two cases.
	\begin{itemize}
		\item[$\bullet$] The components of $\mu$ are either all integers or all half-integers. In this case,  $\mathcal{H}(\mu)=0$. Note that $t_i<t_{i+1}$ if and only if $u_i=u_{i+1}$. We have
		\[\begin{tikzcd}(t_i,t_i+1)\arrow[r,"P_{\C_\mfg}"]& (t_i+\frac{1}{2},t_i+\frac{1}{2}).\end{tikzcd}\]
		The above operation decreases the value of $\Vert \lambda^w(\mu)\Vert^2$ by a total of $\M(\mu)/2$. After these averaging operations, precisely $\caN(\mu)$ entries remain negative and equal $-\frac{1}{2}$.  Then
		\[\begin{tikzcd}-\frac{1}{2}\arrow[r,"P_{\C_\mfg}"]& 0.\end{tikzcd}\]
		The above operation decreases the value of $\Vert \lambda^w(\mu)\Vert^2$ by a total of $\N(\mu)/4$. These two operations map $\lambda^w(\mu)$ to $\lambda_a^w(\mu)$. Thus
		\[\Vert \lambda_a^w(\mu)\Vert^2=\Vert \lambda^w(\mu)\Vert^2-\frac{\M(\mu)}{2}-\frac{\N(\mu)}{4}.\]
		\item[$\bullet$] The first $p$ components of $\mu$ are integers or half-integers, and the remaining $q$ components are of the opposite type. In this case, $\M(\mu)=0$ and $u_i=n-i$ implies that $u_i-u_{i+1}=1/2$ for $1\le i<n$. Note that $t_i<t_{i+1}$ if and only if $u_i=u_{i+1}+1/2$. We have
		\[\begin{tikzcd}(t_i,t_i+\frac{1}{2})\arrow[r,"P_{\C_\mfg}"]& (t_i+\frac{1}{4},t_i+\frac{1}{4})\end{tikzcd}\quad \mbox{and}\quad \begin{tikzcd}(-\frac{1}{2},0)\arrow[r,"P_{\C_\mfg}"]& (0,0).\end{tikzcd}\]
		Moreover, if $u_n=0$, then $t_n=-\frac{1}{2}$. In this case,
		\[\begin{tikzcd}t_n=-\frac{1}{2}\arrow[r,"P_{\C_\mfg}"]& 0.\end{tikzcd}\]
		The above operations decrease the value of $\Vert \lambda^w(\mu)\Vert^2$ by a total of $\caH(\mu)/8$ and $\N(\mu)/4$, respectively. These three operations also map $\lambda^w(\mu)$ to $\lambda_a^w(\mu)$. Thus
		\[\Vert\lambda_a^w(\mu)\Vert^2=\Vert \lambda^w(\mu)\Vert^2-\frac{\N(\mu)}{4}-\frac{\caH(\mu)}{8}.\]
	\end{itemize}
	Combining the above discussion, we obtain \eqref{lambda norm-strat}.
\end{proof}

\begin{cor}\label{co-lambda-L}
Take $\mu$ as in \eqref{mu}. Then
	\begin{equation}\label{Diff-1}
		\Vert \lambda^w(\mu)\Vert^2-L(\mu)\le \min\left\{\frac{p}{2},\frac{2q+1}{4}\right\}.
	\end{equation}
\end{cor}
\begin{proof}
	 	Choose $w\in W(\mfg,\mft)^1$ such that $\mu+2\rho_c\in \C_w$. By Proposition \ref{lambda norm}, the left-hand side of \eqref{Diff-1} equals
	\[\frac{\caM(\mu)}{2}+\frac{\caN(\mu)}{4}+\frac{\caH(\mu)}{8}.\]
	Write $\mu+2\rho_c$ and $w^{-1}(\mu+2\rho_c)$ as $(x_1,\dots,x_p\mid y_1,\dots,y_q)$ and $(u_1,\dots,u_n)$, respectively. Then $(u_1,\dots,u_n)$ is merely the decreasing rearrangement of $(x_1,\dots,x_p\mid y_1,\dots,y_q)$. If $u_l$ is a rearrangement of the $x_i$ (resp., $y_j$), then we say that the \emph{source} of $u_l$ is $x$ (resp., $y$).
	 If the index $1\le l\le n-1$ contributes to $\caM(\mu)$, then $i$ corresponds to a pair $(u_{l}, u_{l+1})$  with $u_{i_1}=u_{i_1+1}$. If the index $1\le l\le n-1$  contributes to $\caH(\mu)$, then $l$ corresponds to a pair $(u_{l}, u_{l+1})$ with $u_{l}=u_{l+1}+\frac{1}{2}$. The pairs indicated above are necessarily pairwise distinct, and each pair consists of one element from $x$ and one element from
	 $y$. Indeed, recall that for distinct index $i$ and $j$,
	 \[|x_i-x_j|\ge 2\quad\mbox{and}\quad |y_i-y_j|\ge 2.\]
	  Then a difference of less than $1$ between $u_l$ and $u_{l+1}$ implies that they come from different sources. Moreover, when two pairs intersect, we have triple $(u_l,u_{l+1},u_{l+2})$. The difference between any two of them is at most $1$. This is impossible, because the three elements cannot all come from different sources. Now if the index $1\le l\le n-1$  contributes to $\caH(\mu)$, then $u_l$ and $u_{l+1}$ have different source, since $u_{l+1}\ge n-l+1$ and $|u_l-u_{l+1}|\le 1$. Furthermore, $l+1$ does \emph{not} contributes to $\caM(\mu)$ and $\caN(\mu)$. Otherwise,
	  \[|u_l-u_{l+2}|\le |u_l-u_{l+1}|+|u_{l+1}-u_{l+2}|\le\frac{3}{2}.\]
	  This shows that $u_l$ and $u_{l+2}$ come from different sources, which contradicts the fact that $u_{l+1}$ and $u_{l+2}$ come from different sources.
	
	  Since each index contributes to at most one of $\caM(\mu)$, $\caN(\mu)$ and $\caH(\mu)$, we have
	  \[\caH(\mu)\le 2\big(p-\caM(\mu)-\caN(\mu)\big).\]
	   Note that only $n$ contributes to $\caH(\mu)$ only when $u_n=0$,  and in this case the source of $u_n$ is $x$.  Then
	   \[\caH(\mu)\le 2\big(q-\caM(\mu)-\caN(\mu)\big)+1.\]
	   Now
	   \[\frac{\caM(\mu)}{2}+\frac{\N(\mu)}{4}+\frac{\caH(\mu)}{8}\le \frac{1}{2}\left(\caM(\mu)+\N(\mu)+\frac{\caH(\mu)}{2}\right)\le \min\left\{\frac{p}{2},\frac{2q+1}{4}\right\},\]
	   as desired.
\end{proof}

Note that $\lambda^w(\mu)$ is the same in the cases $a_p>0$ and $a_p<0$. This shows that $\Vert \mathrm{sgn}(\mu)\Vert_{\mathrm{lambda}}=\Vert \mu\Vert_{\mathrm{lambda}}$.

\subsection{Estimate of the lambda norm}
Let
\begin{equation}\label{kappa}
\kappa=\max\{q-p,p-q-1\}.
\end{equation}
The number $\kappa$ reflects the fact that we have absorbed two classes of simple Lie algebras into one in \eqref{g0}.

For any integer $t$, set
\begin{equation}\label{S}
	\caS(t)=\frac{t(4t^2-1)}{12}.
\end{equation}
When $t\ge 0$,
\[\caS(t)=\sum_{r=1}^t \left(r-\frac{1}{2}\right)^2.\]
\begin{lemma}\label{lambda-lem}
	Take $\mu$ as in \eqref{mu}. Then
	\[L(\mu)=\max_{\eta\in \C_\mfg}\big(2\langle \lambda^w(\mu),\eta\rangle-\langle\eta,\eta\rangle\big).\]
	In particular, if there exists $\eta_0\in \C_\mfg$ such that $\lambda^w(\mu)\gg \eta_0$, then
	\[L(\mu)\ge \langle\eta_0,\eta_0\rangle.\]
\end{lemma}

\begin{proof}
	Note that $C_\mfg$ is a convex cone. Then $\lambda_a^w(\mu)$ is orthogonal to $\lambda^w(\mu)-\lambda_a^w(\mu)$. We have
	\begin{align*}
		L(\mu)&=\langle \lambda_a^w(\mu), \lambda_a^w(\mu)\rangle\\
		&=\langle \lambda^w(\mu), \lambda^w(\mu)\rangle-\langle\lambda^w(\mu)-\lambda_a^w(\mu), \lambda^w(\mu)-\lambda_a^w(\mu)\rangle\\
		&= \max_{\eta\in \C_\mfg}\big(2\langle \lambda^w(\mu),\eta\rangle-\langle\eta,\eta\rangle\big).
	\end{align*}
	If exists $\eta_0\in \C_\mfg$ such that $\lambda^w(\mu)\gg \eta_0$, then $\langle \lambda^w(\mu), \eta_0\rangle\ge \langle\eta_0, \eta_0\rangle$. This shows that $L(\mu)\ge \langle\eta_0,\eta_0\rangle$.
\end{proof}

\begin{cor}\label{lambda-co}
	Take $\mu$ as in \eqref{mu}. Let $u_1$ be the first component of $\lambda^w(\mu)$. Then
	\begin{equation}\label{u1}
		u_1=\max\{a_1+p-q-\frac{3}{2},b_1+q-p-\frac{1}{2}\}.
	\end{equation}
	For any nonnegative number $\ell$ satisfying $\ell\le u_1$,
	\[L(\mu)\ge\ell^2+\caS(\kappa-1).\]
	 In particular,
	 \begin{equation}\label{L-lowest}
	 	L(\mu)\ge \caS(\kappa-1)+\frac{1}{4}.
	 \end{equation}
\end{cor}
\begin{proof}
	Note that \eqref{u1} is obvious. For the next inequality, without loss of generality, we assume that $\ell\ge \kappa-\frac{3}{2}$. Write $\lambda^w(\mu)$ as $(u_1,\dots,u_n)$. Set
	\[\tau_\ell=\left\{
	\begin{array}{ll}
		(\ell,\kappa-\frac{3}{2}, \kappa-\frac{5}{2},\dots,\frac{1}{2},\underbrace{0,\dots,0}_{n-\kappa}),& \mbox{if $\kappa\ge 1$},\\
		(\ell,0,\dots,0), &\mbox{if $\kappa=0$}.
	\end{array}
	\right.\]
	It follows from \eqref{lambda-w-triv} below that for $1\le i\le \kappa$,
	\[u_i\ge \kappa-i+\frac{1}{2}.\]
	Thus $u\gg \tau_\ell$. By Lemma \ref{lambda-lem}, we have
	\[L(\mu)\ge \Vert \tau_\ell\Vert^2=\ell^2+\sum_{i=1}^{\kappa-1}\left(i-\frac{1}{2}\right)^2=\ell^2+\caS(\kappa-1)+\frac{\delta_{\kappa,0}}{4}.\]
	Note that $\kappa\ge 1$ impiles $u_1\ge 1/2$, then
	\[L(\mu)\ge \caS(\kappa-1)+\frac{1}{4}.\]
\end{proof}

\subsection{The number $B_0$ and the trivial $K$-type}  We  the trivial $K$-type by $\texttt{triv}$. It is not hard to compute that
\begin{equation}\label{lambda-w-triv}
\lambda^w(\texttt{triv})=\bigg(\underbrace{\kappa-\frac{1}{2},\kappa-\frac{3}{2},
\kappa-\frac{5}{2},\dots,\frac{1}{2}}_{\kappa},-\frac{1}{2}, \dots, -\frac{1}{2}\bigg).
\end{equation}
and that
\[
\lambda_a^w(\texttt{triv})=\bigg(\underbrace{\kappa-\frac{1}{2},\kappa-\frac{3}{2},\kappa-\frac{5}{2},\dots,\frac{1}{2}}_{\kappa},0,\dots,0\bigg).
\]
Thus
\[L(\texttt{triv})=\sum_{r=1}^{\kappa}\left(r-\frac{1}{2}\right)^2=\caS(\kappa).\]
 Note that
 \begin{equation}\label{B0-S}
 	B_0=\caS(n-1)-\caS(\kappa-1).
\end{equation}
 Then
\begin{equation}\label{triv}
	B_0+L(\texttt{triv})=\caS(n-1)+\left(\kappa-\frac{1}{2}\right)^2.
\end{equation}

Using techniques similar to \cite{DY},  we find that
 \[
 \min_{m\in\bbZ_{\ge 0}}\Vert \chi_0+m\beta\Vert_{\mathrm{spin}}=\left\{
 \begin{array}{ll}
 	\Vert p\beta\Vert_{\mathrm{spin}}, & \mbox{if $q\ge p\ge 2$},\\
 	\Vert q\beta\Vert_{\mathrm{spin}}, & \mbox{if $p=q+1$},\\
 	\Vert (q+1)\beta\Vert_{\mathrm{spin}}, & \mbox{if $p>q+1$}.
 \end{array}
 \right.
 \]
Moreover, all the above spin norms are  attained at $(q-\frac{1}{2},\dots,q-\frac{1}{2}\mid p,0,\dots,0)\in\Omega_{p, q}^+$. Note that  $\min_{m\in\bbZ_{\ge 0}}\Vert \chi_0+m\beta\Vert_{\mathrm{spin}}^2$ is precisely the right hand side of \eqref{triv}. Thus
\begin{equation}\label{B0-explanation}
B_0=B(\texttt{triv}).
\end{equation}
See \eqref{B-mu} for the definition of $B(\mu)$.

\section{Properties of $D_w(\mu)$}\label{sec-Dw-mu}

Previously, we found the convexity of $D_w(\mu)$ for $Sp(p, q)$ in \cite{DY2}. Now let us extend it to the general case.

\begin{prop}\label{DW} Let $G$ be a connected simple real non-compact Lie group with finite center.
	  Let $\mu$ be the highest weight of a $K$-type. Fix a $w\in W(\mfg,\mft)^1$. Then $D_{w}(\mu)$ is convex in $\mu$.
\end{prop}
\begin{proof}
	By Proposition \ref{lambda-rho}, we have
	\begin{equation}\label{Dw-1}
		\Vert\mu+2\rho_c\Vert^2-\Vert\mu\Vert_{\mathrm{lambda}}^2=\max_{\eta\in \C_\rho}(2\langle \mu+2\rho_c,\eta\rangle-\Vert\eta\Vert^2).
	\end{equation}
	Note that
	\begin{equation}\label{Dw-2}
		\langle \{\mu-\rho_n^{(w)}\},\rho_c\rangle =\max_{\widetilde{w}\in W(\mfk,\mft)}\langle \mu-\rho_n^{(w)},\widetilde{w}\rho_c\rangle,
	\end{equation}
	By \eqref{Dw-2}, we have
	\begin{equation}\label{Dw-3}
		\Vert \{\mu-\rho_n^{(w)}\}+\rho_c\Vert^2=\Vert\mu-\rho_n^{(w)}\Vert^2+2\max_{\widetilde{w}\in W(\mfk,\mft)}\langle \mu-\rho_n^{(w)},\widetilde{w}\rho_c\rangle+\Vert \rho_c\Vert^2.
	\end{equation}
	Then \eqref{Dw-1} and \eqref{Dw-3} gives
	\begin{equation}\label{Dw}
		D_w(\mu)=A(\mu)+2\max_{\widetilde{w}\in W(\mfk,\mft)}\langle \mu-\rho_n^{(w)},\widetilde{w}\rho_c\rangle+\max_{\eta\in \C_\rho}(2\langle \mu+2\rho_c,\eta\rangle-\Vert\eta\Vert^2)+\Vert \rho_c\Vert^2,
	\end{equation}
where
$$
A(\mu)=\Vert\mu-\rho_n^{(w)}\Vert^2-\Vert\mu+2\rho_c\Vert^2
$$
is linear in $\mu$. All the  terms in \eqref{Dw} are convex in $\mu$, as desired.
\end{proof}

\begin{rmk}
The above proposition actually holds for any $\Delta^+(\frk, \frt)$-dominant weight.
\end{rmk}

Take $\mu$ as in \eqref{mu}. Choose $w\in \Delta(\mfg,\mft)^1$ such that $\mu+2\rho_c\in \C_w$.
Let us deduce an expression for $D_w(\mu)$ that will be used later.  
 Then
\[\Vert\mu-\rho^{(w)}_n+\rho_c\Vert=\Vert\lambda^w(\mu)\Vert.\]
Thus
\[\Vert \{\mu-\rho^{(w)}_n\}+\rho_c\Vert^2-\Vert\mu-\rho^{(w)}_n+\rho_c\Vert^2=D_w(\mu)+L(\mu)-\Vert\lambda^w(\mu)\Vert^2.\]
On the other hand, let $\{\mu-\rho_n^{(w)}\}=\widetilde{w}(\mu-\rho_n^{(w)})$. Then
\begin{align*}
	\Vert \{\mu-\rho^{(w)}_n\}+\rho_c\Vert^2-\Vert\mu-\rho^{(w)}_n+\rho_c\Vert^2&=\langle\{\mu-\rho_n^{(w)}\}-(\mu-\rho_n^{(w)}),2\rho_c \rangle\\
	&=\sum_{\alpha\in\Delta^+(\mfk,\mft)}\langle\mu-\rho_n^{(w)}, \widetilde{w}\alpha\rangle-\sum_{\alpha\in\Delta^+(\mfk,\mft)}\langle\mu-\rho_n^{(w)}, \alpha\rangle\\
	&=2\sum_{\alpha\in\Delta^+(\mfk,\mft)}\max\{\langle\rho_n^{(w)}-\mu,\alpha\rangle,0\}\\
	&=4\langle\rho_n^{(w)},\rho_c\rangle-2\sum_{\alpha\in\Delta^+(\mfk,\mft)}\min\{\langle\mu,\alpha\rangle,\langle\rho_n^{(w)},\alpha\rangle\}.
\end{align*}
Set
\begin{equation}\label{Fw}
	\F_w(\mu)=\sum_{\alpha\in\Delta^+(\mfk,\mft)}\min\{\langle\mu,\alpha\rangle,\langle\rho_n^{(w)},\alpha\rangle\},
\end{equation}
we obtain that
\begin{equation}\label{Dw-decom}
	D_w(\mu)=4\langle \rho_n^{(w)},\rho_c\rangle -2\F_w(\mu)+\Vert\lambda^w(\mu)\Vert^2-L(\mu).
\end{equation}

\section{Proof of Theorem \ref{thm-main}}\label{sec-proof}
Take $\mu$ as in \eqref{mu}. To prove Theorem \ref{thm-main}, we move on according to the ranges of $a_1$ and $b_1$. In each case, we find a suitable $m\in\bbZ_{\ge 0}$ and $w\in W(\mfg,\mft)^1$ such that
\[D_w(\mu+m\beta)\le B_0,\]
thereby showing that $B(\mu)\le B_0$. Note that
\begin{equation}\label{D}
	D_w(\mu+m\beta)= \Vert\{\mu+m\beta-\rho_n^{(w)}\}+\rho_c \Vert^2-L(\mu).
\end{equation}
Let $u_1$ be the first component of $\lambda^w(\mu)$. We will find a $\Delta^+(\mfk,\mft)$-dominant weight $\eta$ such that $\eta\succeq_\mfk \{\mu+m\beta-\rho_n^{(w)}\}$.
Then by \eqref{compare-spin},
$$
\Phi(\eta)=\Vert \eta+\rho_c \Vert^2 \ge\Vert\{\mu+m\beta-\rho_n^{(w)}\}+\rho_c \Vert^2.
$$
On the other hand, we will find an $\ell$ such that $\ell\le u_1$ by means of \eqref{u1}. Then  Corollary \ref{lambda-co} shows that $L(\mu)\ge \ell^2+\caS(\kappa-1)$.  Thus
\[D_w(\mu+m\beta)\le \Phi(\eta)-L(\mu)\le \Phi(\eta)-\ell^2-\caS(\kappa-1).\]
Thus to prove $D_w(\mu+m\beta)\le B_0$, it suffices  to show that
$$
\Phi(\eta)-\ell^2-\caS(\kappa-1)\le B_0.
$$
In what follows, we will use this directly without further mention, especially the use of $L(\mu)\ge \ell+\caS(\kappa-1)$ for inequality estimates.

\subsection{Comparison of weights}\label{CW}

In this subsection, by finding special weights and elements in $\Omega_{p,q}^+$ and comparing them with $B_0$, we obtain several propositions, which establish Theorem \ref{thm-main} for $a_1\le q$ and $b_1<p$.

For nonnegative integer or half-integer $s$, set
\begin{equation}\label{delta-s}
	\delta_s=\left\{
	\begin{array}{ll}
		1, & \mbox{if $s\in \textbf{I}$},\\
		0, & \mbox{otherwise}.
	\end{array}
	\right.
\end{equation}

\begin{lemma}\label{rho-n-h}
	Take $\mu$ as in \eqref{mu}. Assume in addition that $a_1\le q$ and $b_1< p$. For any given integer $0\le h\le p-1$, there exists $\rho_n^{(w)}:=(K_w\mid R_w)=(k_1,\dots,k_p\mid r_1,\dots,r_q)\in \Omega_{p,q}^+$ such that $r_1$ can be taken to be $p-1$ or $p$, and that for $2\le j \le q$,
	\begin{equation}\label{rho-n-h-1}
		|b_j-r_j|\le\left\{
		\begin{array}{ll}
			h, & \mbox{if $b\in \emph{\textbf{I}}$},\\
			h+\frac{1}{2}, & \mbox{if $b\in \emph{\textbf{H}}$},
		\end{array}
		\right.
	\end{equation}
	and there is a nonnegative number $T\le p-1-h$ such that
	\begin{equation}\label{rho-n-h-2*}
		\sum_{i=2}^p |a_i-k_i|\le (q-1)T+\frac{p-1}{2}\delta_{a_1},
	\end{equation}
	where $\delta_{a_1}$ is defined in \eqref{delta-s}. In particular,
	\begin{equation}\label{rho-n-h-2}
		\sum_{i=2}^p |a_i-k_i|\le (q-1)(p-1-h)+\frac{p-1}{2}\delta_{a_1}.
	\end{equation}
\end{lemma}
\begin{proof}
	Set $\widetilde{a}_i=\max\{0,\lceil a_i\rceil-1\}$ for $2\le i\le p$ and $\widetilde{b}_j=\min\{p-1,\lceil b_j\rceil\}$ for $2\le j\le q$, where $\lceil\cdot\rceil$ denotes the ceiling function. Let $(q,\widetilde{a}_2\dots,\widetilde{a}_p\mid y_1,\dots,y_q)$ be the unique element of $\Omega_{p,q}^*$ determined by $(q,\widetilde{a}_2,\dots,\widetilde{a}_p)$. For $2\le j\le q$, take
	\[r_j=\left\{
	\begin{array}{ll}
		\widetilde{b}_j-h, & \mbox{if $y_j<\widetilde{b}_j-h$},\\
		y_j, & \mbox{if $\widetilde{b}_j-h\le y_j\le\widetilde{b}_j+h$},\\
		\widetilde{b}_j+h, & \mbox{if $\widetilde{b}_j+h< y_j$}.
	\end{array}
	\right.\]
	Let $r_1$ be $p-1$ or $p$. Then $R_w=(r_1,\dots,r_q)$  determines a unique element $\rho_n^{(w)}$ in $\Omega_{p,q}^+$. Then \eqref{rho-n-h-1} holds immediately. Note that $|q+\frac{1}{2}-k_1|=|y_1-r_1|$. Define
	\[T=\max_{2\le j\le q}|y_j-r_j|.\]
	Then $T\le p-1
	-h$ and Lemma \ref{Omega-pq} shows that
	\begin{equation}\label{iq-1}
		\sum_{i=2}^p |\widetilde{a}_i+\frac{1}{2}-k_i|=\sum_{j=2}^q|y_j-r_j|\le (q-1)T\le (q-1)(p-1-h).
	\end{equation}
	On the other hand, we have
	\[|a_i-k_i|\le |\widetilde{a}_i+\frac{1}{2}-k_i|+\frac{\delta_{a_1}}{2}\]
	for $2\le i\le p$. Summing the above inequalities from $i=2$ to $p$ and using \eqref{iq-1}, we obtain \eqref{rho-n-h-2*} and \eqref{rho-n-h-2}.
\end{proof}
For any $\mu$ in \eqref{mu} satisfying $a_1\le q$ and $b_1<p$, and for any $0\le h\le p-1$, Lemma \ref{rho-n-h} yields a finite subset of $\Omega_{p,q}^+$ satisfying \eqref{rho-n-h-1} and \eqref{rho-n-h-2*}. For convenience, we denote this finite subset by $\Omega_{\mu}(h)$.

\begin{prop}\label{pro-1}
	Take $\mu$ as in \eqref{mu}. Assume that $a_1<q$ and $b_1<p$. If $\lfloor q-a_1\rfloor\ge \lfloor p-b_1\rfloor$, then $B(\mu)\le B_0$. Here $\lfloor \cdot\rfloor$ denotes the floor function.
\end{prop}
\begin{proof}
	By Lemma \ref{rho-n-h}, take $\rho_n^{(w)}:=(k_1,\dots,k_p\mid r_1,\dots,r_q)\in \Omega_{\mu}(0)$ such that $r_1=p$. We consider the following four cases.
	
	(i) $\mu\in (\textbf{I}\mid\textbf{I})$. If $\lfloor q-a_1\rfloor=\lfloor p-b_1\rfloor$, we have
	\[\{\mu+\lfloor p-b_1\rfloor\beta-\rho_n^{(w)}\}\preceq_\mfk \eta_1 :=\left.\left(q-\frac{1}{2},\dots,q-\frac{1}{2},\frac{1}{2}\ \right|\ 0,\dots,0\right).\]
	 By \eqref{L-lowest}, we have $L(\mu)\ge \caS(\kappa-1)+1/4$. Then
	\[\Phi(\eta_1)-L(\mu)=\caS(n-1)-\frac{1}{4}-L(\mu)=\caS(n-1)-S(\kappa-1)-\frac{1}{4}=B_0-\frac{1}{4},\]
	where the last step uses \eqref{B0-S}. If $\lfloor q-a_1\rfloor>\lfloor p-b_1\rfloor$, we have	
		\[\{\mu+\lfloor p-b_1\rfloor\beta-\rho_n^{(w)}\}\preceq_\mfk \eta_1^* :=\left.\left(q-\frac{1}{2},\dots,q-\frac{1}{2},q-\frac{1}{2}-\lfloor p-b_1\rfloor\ \right|\ 0,\dots,0\right).\]
	Take $\ell=q-\frac{1}{2}-\lfloor p-b_1\rfloor$. Then
	\[\Phi(\eta_1^*)-L(\mu)=\caS(n-1)+\ell^2-L(\mu)\le \caS(n-1)-\caS(\kappa-1)= B_0.\]
	
	(ii) $\mu\in (\textbf{H}\mid \textbf{I})$. We have
	\[\{\mu+\lfloor p-b_1\rfloor\beta-\rho_n^{(w)}\}\preceq_\mfk \eta_2 :=\left.\left(q-1,\dots,q-1,q-1-\lfloor p-b_1\rfloor\ \right|\ 0,\dots,0\right).\]
	Take $\ell=q-1-\lfloor p-b_1\rfloor$. Then
	\begin{align*}
		\Phi(\eta_2)-L(\mu)&=\caS(n-1)+\ell^2-(p-1)\left(q+\frac{p}{2}-\frac{3}{4}\right)-L(\mu)\\
		&\le B_0-(p-1)\left(q+\frac{p}{2}-\frac{3}{4}\right)\le B_0.
	\end{align*}
	
	(iii) $\mu\in (\textbf{I}\mid \textbf{H})$. In this case, we \emph{claim} that
	\[\sum_{i=2}^p|a_i-k_i|\le (p-2)\left(q-\frac{1}{2}\right)+\max\{a_1-1,0\}+\frac{1}{2}.\]
	Indeed, in the notation of Lemma \ref{rho-n-h}, we have $y_j=p-1$ for $1\le j\le q-a_1+1$. Then $\widetilde{b}_j=\min\{p-1,\lceil b_j\rceil\}\ge 1$ shows that
	\[|y_j-r_j|=|y_j-\widetilde{b}_j|\le p-2\]
	for $1\le j\le q-a_1+1$. Therefore,
	\begin{align*}
		\sum_{i=2}^p|a_i-k_i|&\le \sum_{j=2}^q|y_j-r_j|+\frac{p-1}{2}\\
		&\le (p-2)(q-a_1+1)+(p-1)\max\{a_1-1,0\}+\frac{p-1}{2}\\
		&=\left(p-2\right)\left(q-\frac{1}{2}\right)+\max\{a_1-1,0\}+\frac{1}{2}.
	\end{align*}
	Thus the claim holds.
	
	Note that $k_1=q-\frac{1}{2}$. If $\lfloor p-b_1\rfloor=0$, we have
	\[\{\mu+\beta-\rho_n^{(w)}\}\preceq_\mfk \eta_3:=\left.\left(q-\frac{1}{2},\dots,q-\frac{1}{2},q-\frac{3}{2},\frac{1}{2}\ \right|\ \frac{1}{2},\dots,\frac{1}{2}\right).\]
	Take $\ell=q-1$. Then
	\begin{align*}
		\Phi(\eta_3)-L(\mu)&=\caS(n-1)+\frac{2q^2-7q+1}{4}-L(\mu)\\
		&\le B_0+\frac{2q^2-7q+1}{4}-\ell^2\\
		&= B_0-\frac{2q^2-q+3}{4}\le B_0.
	\end{align*}
	If $\lfloor p-b_1\rfloor>0$, we have
	\[\{\mu+\lfloor p-b_1\rfloor\beta-\rho_n^{(w)}\}\preceq_\mfk \eta_3^* :=\left.\left(q-\frac{1}{2},\dots,q-\frac{1}{2},q-\frac{1}{2}-\lfloor p-b_1\rfloor,\frac{1}{2}\ \right|\ \frac{1}{2},\dots,\frac{1}{2}\right).\]
	Take $\ell=q-1-\lfloor p-b_1\rfloor$. Then
	\begin{align*}
		\Phi(\eta_3^*)-L(\mu)&=\caS(n-1)+\lfloor p-b_1\rfloor^2-(2q+1)\lfloor p-b_1\rfloor+\frac{2q^2+q+1}{4}-L(\mu)\\
		&\le B_0+\frac{2q^2+q+1}{4}+\lfloor p-b_1\rfloor^2-\ell^2-(2q+1)\lfloor p-b_1\rfloor\\
		&=B_0-3\lfloor p-b_1\rfloor-\frac{2q^2-9q+3}{4}\le B_0.
	\end{align*}
	
	(iv) $\mu\in (\textbf{H}\mid \textbf{H})$. Since $r_2<p$ and $r_j\ge 1$ for $1\le j\le q$, then $k_1=q-\frac{1}{2}$ and $k_p=\frac{1}{2}$. Note that $|a_1+\lfloor p-b_1\rfloor-k_1|=\lfloor q-a_1\rfloor-\lfloor p-b_1\rfloor$ and $|a_p-k_p|\le q-1-\lfloor q-a_1\rfloor$. We have
	\[\{\mu+\lfloor p-b_1\rfloor\beta-\rho_n^{(w)}\}\preceq_\mfk \eta_4 :=\bigg(q-1,\dots,q-1,q-1-\lfloor p-b_1\rfloor,0\ \left|\ \frac{1}{2},\dots,\frac{1}{2}\bigg)\right..\]
	Take $\ell=q-1-\lfloor p-b_1\rfloor$. Then
	\begin{align*}
		\Phi(\eta_4)-L(\mu)&=\caS(n-1)+(\ell+1)^2-\frac{(n-1)(2n-3)}{4}-L(\mu)\\
		&\le B_0+(\ell+1)^2-\ell^2-\frac{(n-1)(2n-3)}{4}\\
		&=B_0+\ell(2\ell+1)-\frac{(n-1)(2n-3)}{4}\\
		&\le B_0-\frac{2q^2-5q+5}{4}-2\lfloor p-b_1\rfloor\le B_0,
		\end{align*}
where the last step uses $n=p+q\ge q+2$.

This finishes the proof.
\end{proof}

\begin{prop}\label{pro-2}
	Take $\mu$ as in \eqref{mu}. Assume that $a_1<q$ and $b_1<p$. If $\lfloor q-a_1\rfloor<\lfloor p-b_1\rfloor$, then $B(\mu)\le B_0$.
\end{prop}
\begin{proof}
	Write $\lfloor p-b_1\rfloor-\lfloor q-a_1\rfloor$ as $s$. Let $h_1=s-2+\delta_{1,s}$, $h_2=s-1$ and $h_3=p-1$. By Lemma \ref{rho-n-h}, we can take $\rho_n^{(w_1)}\in \Omega_{\mu}(h_1)$, $\rho_n^{(w_2)}\in \Omega_{\mu}(h_2)$ and $\rho_n^{(w_3)}\in \Omega_{\mu}(h_3)$ such that their $(p+1)$-th components are $p-1$. We proceed according to  the following four cases.

	(i) $\mu\in (\textbf{I}\mid \textbf{I})$. In this case,  from \eqref{rho-n-h-2},
		\[\{\mu+\big(\lfloor q-a_1\rfloor+1-\delta_{1,s}\big)\beta-\rho_n^{(w_{1})}\}\preceq_\mfk\eta_1:=\bigg(\underbrace{q-\frac{1}{2},\dots,q-\frac{1}{2}}_{p-1-h_1},\dfrac{1}{2},\dots,\dfrac{1}{2}\ \bigg|\ h_1,\dots,h_1\bigg).\]
		If $h_1>0$, then $s=h_1+2$. Thus
		\[a_1+p-q-\frac{1}{2}\ge b_1+h_1+\frac{1}{2}\ge h_1+\frac{1}{2}.\]
		Take $\ell=h_1+\frac{1}{2}$. We have
		\[\Phi(\eta_1)-L(\mu)=\caS(n-1)+\left(h_1+\frac{1}{2}\right)^2-L(\mu)\le B_0.\]
		If $h_1=0$, by \eqref{L-lowest}, we have
		\[\Phi(\eta_1)-L(\mu)=\caS(n-1)+\frac{1}{4}-L(\mu)\le B_0.\]

		(ii) $\mu\in(\textbf{H}\mid \textbf{I})$. In this case, from \eqref{rho-n-h-2},
		\[\{\mu+\big(\lfloor q-a_1\rfloor+1-\delta_{1,s}\big)\beta-\rho_n^{(w_{1})}\}\preceq_\mfk\eta_2=(\underbrace{q-1,\dots,q-1}_{p-1-h_1},\delta_{s,1},0,\dots,0\mid h_1,\dots,h_1).\]
		Similarly, if $h_1>0$, then
			\[a_1+p-q-\frac{1}{2}\ge b_1+h_1\ge h_1.\]
		Now whether or not $h_1>0$, we take $\ell=h_1$. Then
		\begin{align*}
			\Phi(\eta_2)-L(\mu)&= \caS(n-1)+h_1^2-q(p-1-h_1)-\frac{(p-1)(2q-3)}{4}+\delta_{1,s}-L(\mu)\\
			&\le B_0-\left(q(p-1-h_1)+\frac{(p-1)(2q-3)}{4}-\delta_{1,s}\right)\le B_0,
		\end{align*}
		where the last step uses $0\le h_1\le p-2$.
		
		(iii) $\mu\in (\textbf{I}\mid \textbf{H})$. In this case, we \emph{claim} that $T$ in \eqref{rho-n-h-2*} is at most $p-2-h_2$ for $\Omega_{\mu}(h_2)$. Indeed, in the notation of Lemma \ref{rho-n-h}, note that $1\le \widetilde{b}_j\le p-\lfloor p-b_1\rfloor\le p-2$ for $2\le j\le q$ and $y_j\le p-1$. Then
		\[T\le \max_{2\le j\le q}\{|y_j-\widetilde{b}_j|\}-h_2\le p-2-h_2.\]
		Thus the cliam holds. By \eqref{rho-n-h-2*},
		\[\{\mu+(\lfloor q-a_1\rfloor+1)\beta-\rho_n^{(w_{2})}\}\preceq_\mfk\eta_3:=\bigg(\underbrace{q-\frac{1}{2},\dots,q-\frac{1}{2}}_{p-2-h_2},\dfrac{1}{2},\dots,\dfrac{1}{2}\ \bigg|\ h_2+\dfrac{1}{2},\dots,h_2+\dfrac{1}{2}\bigg).\]
		Note that
		\[a_1+p-q-\frac{1}{2}\ge b_1+h_2\ge h_2+\frac{1}{2}.\]
		Take $\ell=h_2+\frac{1}{2}$. Then
		\begin{align*}
			\Phi(\eta_3)-L(\mu)&= \caS(n-1)+h_2^2-(q-3)h_2-\frac{(q+3)(2q-3)}{4}-L(\mu)\\
			&\le B_0-(q-2)h_2-\frac{2q^2+3q-8}{4}\le B_0.
		\end{align*}
		
		(iv) $\mu\in (\textbf{H}\mid \textbf{H})$. In this case,
		since $b_1\le p-\frac{3}{2}$, we have
		\[|b_j-r_j|\le \min\{h_3+1/2, p-3/2\}\]
		for $j\ge 2$.  It is easy to see that $r_j=y_j$ for $j\ge 2$ in the notation of Lemma \ref{rho-n-h}. Then
\[\{\mu+(\lfloor q-a_1\rfloor+1)\beta-\rho_n^{(w_3)}\}\preceq_\mfk \eta_4:=\bigg(0,\dots,0\ \bigg|\  p-\frac{3}{2},\dots,p-\frac{3}{2}\bigg).\]
		
		Since $\lfloor q-a_1\rfloor\le \lfloor p-b_1\rfloor-1\le p-2$, we have
		\[a_1+p-q-\frac{3}{2}=p-\lfloor q-a_1\rfloor-2\ge 0.\]
		Take $\ell=p-\lfloor q-a_1\rfloor-2$. Then
		\begin{align*}
			\Phi(\eta_4)-L(\mu)&=\caS(n-1)+(p-1)^2-\frac{(n-1)(2n-3)}{4}-L(\mu)\\
			&\le B_0-\left(\frac{(n-1)(2n-3)}{4}+\ell^2-(p-1)^2\right).
		\end{align*}
		We compute that
		\begin{align*}
			&\ \quad\frac{(n-1)(2n-3)}{4}+\ell^2-(p-1)^2\\
			&\ge \frac{(\ell+2\lfloor q-a_1\rfloor+2)(2\ell+4\lfloor q-a_1\rfloor+3)}{4}+\ell^2-(\ell+\lfloor q-a_1\rfloor+1)^2\\
			&=\frac{\ell(2\ell-1)}{4}+\frac{(\lfloor q-a_1\rfloor+1)(2\lfloor q-a_1\rfloor+1)}{2}>0,
		\end{align*}
		where the first step uses $n-1\ge p+\lfloor q-a_1\rfloor=\ell+2\lfloor q-a_1\rfloor+2$.

This finishes the proof.
\end{proof}

\begin{prop}\label{pro-3}
	Take $\mu$ as in \eqref{mu}. Assume that $a_1=q$ and $b_1<p$. Then $B(\mu)\le B_0$.
\end{prop}
\begin{proof}
 We move on according to the following two cases.
	
	(i) $\mu\in (\textbf{I}\mid \textbf{I})$.
	Write $p-b_1-2+\delta_{1,p-b_1}$ as $h_1$. Then $0\le h_1\le p-2$. Take $\rho_n^{(w_1)}\in \Omega_{\mu}(h_1)$ such that its $(p+1)$-th component is $p-1$. By \eqref{rho-n-h-2}, we have
		\[\{\mu+(1-\delta_{1,p-b_1})\beta-\rho_n^{(w_1)}\}\preceq_\mfk \eta_1:=\bigg(\underbrace{q-\frac{1}{2},\dots,q-\frac{1}{2}}_{p-1-h_1},\frac{1}{2},\dots,\frac{1}{2}\ \bigg|\ h_1,\dots,h_1\bigg).\]
			Note that
		\[a_1+p-q-\frac{3}{2}=p-\frac{3}{2}.\]	
		Take $\ell=p-\frac{3}{2}$, we obtain $L(\mu)\ge \left(p-\frac{3}{2}\right)^2+\caS(\kappa-1)$.	Then
		\[\Phi(\eta_1)-L(\mu)\le B_0+\left(h_1+\frac{1}{2}\right)^2-\left(p-\frac{3}{2}\right)^2\le B_0.\]
		
	(ii) $\mu\in (\textbf{I}\mid \textbf{H})$. Firstly, assume that $b_1\le p-\frac{3}{2}$. Let $K_{w_2}$ be $(a_1+\frac{1}{2},\dots,a_p+\frac{1}{2})$. Then $K_{w_2}$ uniquely determines the element $\rho_n^{(w_2)}:=(K_{w_2}\mid R_{w_2})$ in $\Omega_{p,q}^+$. We have
		\[\{\mu-\rho_n^{(w_2)}\}\preceq_\mfk \eta_2:=\left(\frac{1}{2},\dots,\frac{1}{2}\ \left|\  p-\frac{3}{2},\dots, p-\frac{3}{2}\right)\right..\]
	Then
	\[\Phi(\eta_2)-L(\mu)\le  B_0-(q-2)\left(p+\frac{2q-1}{4}\right)-\frac{3}{2}\le B_0.\]
		
	Now assume that $b_1=p-\frac{1}{2}$. Take $\rho_n^{(w_3)}\in \Omega_{\mu}(0)$. By Lemma \ref{rho-n-h}, we have
		\[\{\mu-\rho_n^{(w_3)}\}\preceq_\mfk \eta_3 :=\left.\left(q-\frac{1}{2},\dots,q-\frac{1}{2},\frac{1}{2}\ \right|\ \frac{1}{2},\dots,\frac{1}{2}\right).\]
	Note that $b_1+q-p-\frac{1}{2}=q-1$. Take $\ell=\max\{p-\frac{3}{2},q-1\}$, we obtain that $L(\mu)\ge \ell^2+\caS(\kappa-1)$.  If $q\ge 5$, then
		\[\Phi(\eta_2)-L(\mu)\le B_0-\left((q-1)^2-\frac{2q^2+q+1}{4}\right)=B_0-\frac{(q-5)(2q+1)+8}{4}\le B_0.\]
	 If $q\le 4$ and $p\ge 5$, then
	 \[\Phi(\eta_2)-L(\mu)\le B_0-\left(\left(p-\frac{3}{2}\right)^2-\frac{2q^2+q+1}{4}\right)\le B_0-\left(\frac{49}{4}-\frac{37}{4}\right)\le B_0.\]
	
	 It remains only to consider the case $2\le p,q\le 4$. In this case, we have $\kappa\le 2$ and $\caS(\kappa-1)\le\frac{1}{4}$. Note that the first and $(p+1)$-th components of $\mu+2\rho_c$ are $n+p-2$ and $n+q-\frac{3}{2}$, respectively. Then the sum of the first two components $u_1$ and $u_2$ of $\lambda^w(\mu)$ is at least
		\[(n+p-2)+\left(n+q-\frac{3}{2}\right)-\left(n-\frac{1}{2}\right)-\left(n-\frac{3}{2}\right)=n-\frac{3}{2}.\]
	Thus
		\[L(\mu)\ge \frac{(u_1+u_2)^2}{2}=\frac{1}{2}\left(n-\frac{3}{2}\right)^2.\]
	We have
	\[\Phi(\eta_2)-L(\mu)\le B_0-\left(\frac{1}{2}\left(n-\frac{3}{2}\right)^2-\frac{1}{4}-\frac{2q^2+q+1}{4}\right)\le B_0-\frac{2q-3}{8}\le B_0,\]
	where the second step uses $n-\frac{3}{2}\ge q+\frac{1}{2}$.
	
This finishes the proof.
\end{proof}

\subsection{Truncation of weights}
In this subsection, we prove Theorem \ref{thm-main} for $a_1\le q$ or $b_1<p-\frac{1}{2}$.

	Take $\mu$ as in \eqref{mu}. Let $a^*=(a_1^*,\dots, a_p^*)$, where
\[a_i^*=\left\{
\begin{array}{ll}
	\min\{a_i, q+1\}, &\mbox{if $a\in\textbf{I}$},\\
	\min\{a_i, q+\frac{1}{2}\}, & \mbox{if $a\in\textbf{H}$}.
\end{array}
\right.\]
Set
\[k_i=\left\{
\begin{array}{ll}
	|a_i^*-\frac{1}{2}|, & \mbox{if $a\in\textbf{I}$},\\
	a_i^*, & \mbox{if $a\in\textbf{H}$}
\end{array}
\right.\]
Then $(k_1,\dots,k_p)$ determines a unique element in $\Omega_{p,q}^+$, which we call \emph{the weight corresponding to} $\mu$ \emph{on the left}, denoted by $\Omega_{\lhd}(\mu)$. Similarly, let $b^*=(b_1,\dots,b_q)$, where
\[b_j^*=\left\{
\begin{array}{ll}
	\min\{b_j, p\}, & \mbox{if $b\in\textbf{I}$},\\
	\min\{b_j, p+\frac{1}{2}\}, & \mbox{if $b\in\textbf{H}$}
\end{array}
\right.\]
Set
\[r_j=\left\{
\begin{array}{ll}
	b_j^*, & \mbox{if $b\in\textbf{I}$},\\
	\min\{p,b_j^*+\frac{1}{2}\}, & \mbox{if $b\in\textbf{H}$}
\end{array}
\right.\]
Then $(r_1,\dots,r_p)$ determines a unique element in $\Omega_{p,q}^+$, which we call \emph{the weight corresponding to} $\mu$ \emph{on the right}, denoted by $\Omega_{\rhd}(\mu)$.

\begin{lemma}\label{JD}
	Take $\mu$ as in \eqref{mu}. Then
	\begin{itemize}
		\item[(i)] For $\rho_n^{(w)}=\Omega_{\rhd}(\mu)$, we have
		\begin{equation}\label{Dw-b}
			D_w(a\mid b)\le D_w(a\mid b^*).
		\end{equation}
		\item[(ii)]  For $\rho_n^{(w)}=\Omega_{\lhd}(\mu)$, we have
		\begin{equation}\label{Dw-a}
			D_w(a\mid b)\le D_w(a^*\mid b).
		\end{equation}
	\end{itemize}
\end{lemma}
\begin{proof}
	Write $\rho_n^{(w)}$ as $(k_1,\dots,k_p\mid r_1,\dots,r_q)$.
	
(i)  Let
	\[\delta=(a\mid b)-(a\mid b^*)=(0,\dots,0\mid b_1-b_1^*,\dots,b_q-b_q^*).\]
	Note that $b_1-b_1^*\ge \cdots\ge b_q-b_q^*$.  By \eqref{Dw},
	\begin{align*}
		D_w(a\mid b)-D_w(a\mid b^*)&\le -2\langle\delta,2\rho_c+\rho_n^{(w)}\rangle+2\max_{\widetilde{w}\in W(\mfk,\mft)}\langle \delta,\widetilde{w}\rho_c\rangle +2\max_{\eta\in\C_\rho}\langle \delta,\eta\rangle\\
		&=-2\langle\delta,2\rho_c+\rho_n^{(w)}\rangle+ 2\langle \delta,\rho_c\rangle +2\sum_{j=1}^q\left(n-i+\frac{1}{2}\right)(b_j-b_j^*)\\
		&=2\sum_{j=1}^q(p-r_j)(b_j-b_j^*)=0.
	\end{align*}
	The last step uses $b_j-b_j^*>0$ only if $p=r_j$ for $1\le j\le q$.
	
(ii) Let
\[\delta=(a\mid b)-(a^*\mid b)=(a_1-a_1^*,\dots,a_p-a_p^*\mid 0,\dots, 0).\]
Note that $a_1-a_1^*\ge\cdots\ge a_p-a_p^*$. By \eqref{Dw},
	\begin{align*}
	D_w(a\mid b)-D_w(a^*\mid b)&\le -2\langle\delta,2\rho_c+\rho_n^{(w)}\rangle+2\max_{\widetilde{w}\in W(\mfk,\mft)}\langle \delta,\widetilde{w}\rho_c\rangle +2\max_{\eta\in\C_\rho}\langle \delta,\eta\rangle\\
	&=-2\langle\delta,2\rho_c+\rho_n^{(w)}\rangle+ 2\langle \delta,\rho_c\rangle +2\langle \delta,\rho\rangle\\
	&= 2\langle\delta, \rho_n-\rho_n^{(w)} \rangle=2\sum_{i=1}^p\left(q+\frac{1}{2}-k_i\right)(a_i-a_i^*)=0.
\end{align*}
The last step uses $a_i-a_i^*>0$ only if $q+\frac{1}{2}=k_i$ for $1\le i\le p$.
\end{proof}

\begin{prop}
	Take $\mu$ as in \eqref{mu}. If $b_1< p-\frac{1}{2}$, then $B(\mu)\le B_0$.
\end{prop}

\begin{proof}
	By Propositions \ref{pro-1}, \ref{pro-2} and \ref{pro-3}, it suffices to consider $a_1>q$ and $b_1<p-\frac{1}{2}$. Take $\rho_n^{(w)}=\Omega_{\lhd}(\mu)$. Then Lemma \ref{JD} shows that
	\[D_w(a\mid b)\le D_w(a^*\mid b).\]
	Therefore, without loss of generality, we assume that $\mu=(a^*\mid b)$.  Recall by \eqref{delta-s},
	\[\delta_a=\left\{
	\begin{array}{ll}
		1,& \mbox{if $a\in \textbf{I}$},\\
		0, & \mbox{if $a\in \textbf{H}$}.
	\end{array}
	\right.\]
	Then
	\[\{\mu-\rho_n^{(w)}\}\preceq_\mfk \eta :=\left(\frac{\delta_a}{2},\dots,\frac{\delta_a}{2}\ \bigg|\ p-1,\dots,p-1\right).\]
	Note that $a_1^*=q+\frac{\delta_a}{2}+\frac{1}{2}$. Then
	\[a_1^*+p-q-\frac{3}{2}=p-1+\frac{\delta_a}{2}.\]
	By \eqref{u1}, take $\ell=p-1+\frac{\delta_a}{2}$. We have
	\begin{align*}
	\Phi(\eta)-L(\mu)&=\caS(n-1)+\frac{(2p-1)(p-1+p\delta_a)}{4}-L(\mu)\\
	&\le B_0+\left(\frac{(2p-1)(p-1+p\delta_a)}{4}-\left(p+\frac{\delta_a}{2}-1\right)^2\right).
	\end{align*}
	A straightforward computation shows that
	\[\frac{(2p-1)(p-1+p\delta_a)}{4}-\left(p+\frac{\delta_a}{2}-1\right)^2\le 0\]
	for $\delta_a\in\{0,1\}$. This finish the proof.
\end{proof}

\begin{prop}
	Take $\mu$ as in \eqref{mu}. If $a_1\le q$, then $B(\mu)\le B_0$.
\end{prop}

\begin{proof}
	By Propositions \ref{pro-1} and \ref{pro-2}, it suffices to consider $a_1\le q$ and $b_1\ge p$. Take $\rho_n^{(w)}=\Omega_{\rhd}(\mu)$. Then Lemma \ref{JD} shows that
	\[D_w(\mu)\le D_w(a\mid b^*).\]
	Therefore, without loss of generality, we assume that $\mu=(a\mid b^*)$. Then
	\[\{\mu-\rho_n^{(w)}\}\preceq_\mfk \eta :=\left(q-\dfrac{1}{2},\dots,q-\dfrac{1}{2}\ \bigg|\ \frac{1-\delta_{b_1}}{2},\dots,\frac{1-\delta_{b_1}}{2}\right).\]
	
	If $b\in \textbf{I}$, then $b_1^*=p$. We have $b_1^*+q-p-\frac{1}{2}=q-\frac{1}{2}$. Take $\ell=q-\frac{1}{2}$. Then
	\[\Phi(\eta)-L(\mu)=\caS(n-1)+\left(q-\frac{1}{2}\right)^2-L(\mu)\le B_0.\]
	
	Now assume that $b\in \textbf{H}$. In this case, $b_1^*=p+\frac{1}{2}$ and $r_j\ge 1$ for $1\le j\le q$. Then $k_p=\frac{1}{2}$. We \emph{claim} that it suffices to consider the cases $a_p=0$ and $a_p=q$. Indeed, if $a_p<q$, let
	\[\widetilde{a}=\frac{q}{q-a_p}\left(a_1-a_p, a_2-a_p, \dots, a_{p-1}-a_p, 0\right)\]
	and
	\[\mu_0=(\widetilde{a}\mid b),\quad \mu_1=(q,\dots,q\mid b).\]
	Then $\mu=(1-\frac{a_p}{q})\mu_0+\frac{a_p}{q}\mu_1$. By Proposition \ref{DW} and its remark,
	\[D_w(\mu)\le \left(1-\frac{a_p}{q}\right)D_w(\mu_0)+\frac{a_p}{q}D_w(\mu_1).\]
	This proves the claim. Now we prove these two cases separately.
	
	(i) $a_p=0$. Note that $|a_p-k_p|=\frac{1}{2}$, and $k_1\le q-\frac{1}{2}$ since $r_1=p$. Then
	\[\{\mu-\rho_n^{(w)}\}\preceq_\mfk\eta^* :=\left(q-\frac{1}{2},\dots,q-\frac{1}{2},\frac{1}{2}\ \bigg|\  \frac{1}{2},\dots,\frac{1}{2}\right).\]
	Since $b_1^*=p+\frac{1}{2}$, then $b_1^*+q-p-\frac{1}{2}=q$. Take $\ell=q$, we have
	\[\Phi(\eta^*)-L(\mu)=\caS(n-1)+\frac{2q^2+q+1}{4}-L(\mu)\le B_0-\frac{2q^2-q-1}{4}\le B_0.\]
	
	(ii) $a_p=q$. Then $a=(q,\dots,q)$. Note that
	\[\Phi(\eta)=B_0+\frac{3}{2}q^2+\caS(\kappa-1)-\frac{3q-1}{4}.\]
	Thus we only need to show that
	\[L(\mu)\ge \frac{3}{2}q^2+\caS(\kappa-1).\]
	
	Let $\chi=(\chi_1,\dots,\chi_n)$, where
	\[\chi_1=\max\{q,\kappa\},\quad \chi_i=\max\{\kappa-i+1,0\}\]
	for $2\le i\le n$. Set
	\[\widetilde{\chi}=\chi+\textbf{1}_n,\]
	where $\textbf{1}_n$ denotes the all-ones $n$-vector. Then $\widetilde{\chi}\in \C_\mfg$.  Write $\lambda^w(\mu)$ as $(u_1,\dots,u_n)$. Note that each component of  $\mu$ is at least $1/2$. Then it follows from \eqref{lambda-w-triv} that $\lambda^w(\mu)\gg \chi$. Now
	\begin{equation}\label{ui}
	\begin{split}
		\sum_{i=1}^n u_i&= \sum_{i=1}^p a_i+\sum_{j=1}^qb_j^*+p(p-1)+q^2-\frac{n^2}{2}\\
		&\ge pq+p+\frac{q}{2}+p(p-1)+q^2-\frac{n^2}{2}\\
		&= \frac{p^2+q^2+q}{2},
	\end{split}
\end{equation}
where the second step uses $b_1^*=p+\frac{1}{2}$ and $b_j^*\ge \frac{1}{2}$ for $j\ge 2$.
	On the other hand,
	\begin{equation}\label{chi}
	\begin{split}
		\sum_{i=1}^n\chi_i(\chi_i-2)&=\chi_1(\chi_1-2)+\sum_{r=1}^{\kappa-1}r(r-2)\\
		&\ge q(q-2)+\sum_{r=1}^{\kappa-1}r(r-2)\\
		&=q(q-2)+\caS(\kappa-1)-\frac{\kappa^2}{2}+\frac{\kappa+1}{4}\\
		&\ge q(q-2)+\caS(\kappa-1)-\frac{(p-q)^2}{2},
	\end{split}
	\end{equation}
	where the second step uses $\chi_1\ge q\ge 2$ and the last step uses $\kappa\le |p-q|$.
	By Lemma \ref{lambda-lem},
	\begin{align*}
		L(\mu)&\ge 2\langle\lambda^w(\mu),\widetilde{\chi}\rangle-\langle\widetilde{\chi},\widetilde{\chi}\rangle\\
		&= 2\sum_{i=1}^n u_i-n+\sum_{i=1}^n\chi_i(\chi_i-2)+2\sum_{i=1}^n\chi_i(u_i-\chi_i)\\
		&\ge 2\sum_{i=1}^n u_i-n+\sum_{i=1}^n\chi_i(\chi_i-2)\\
		&\ge \caS(\kappa-1)+\frac{3}{2}q^2+\frac{(p-2)(p+2q)}{2}\\
		&\ge \caS(\kappa-1)+\frac{3}{2}q^2,
	\end{align*}
	where the third step uses $\lambda^w(\mu)\gg \chi$, and the fourth step uses \eqref{ui} and \eqref{chi}. This finishes the proof.
\end{proof}

\subsection{Remaining cases of Theorem \ref{thm-main}}

In this subsection, we assume that $a_1>q$ and $b_1\ge p-\frac{1}{2}$. Take $\mu$ as in \eqref{mu}. Choose $w\in \Delta(\mfg, \mft)^1$ such that $\mu+2\rho_c\in \C_w$. Then by \eqref{Dw-decom} and Corollary \ref{co-lambda-L},
\begin{equation}\label{Dw1}
	D_w(\mu)\le 4\langle \rho_n^{(w)},\rho_c\rangle-2\F_w(\mu)+\min\left\{\frac{p}{2}, \frac{2q+1}{4}\right\}.
\end{equation}

Write $\rho_n^{(w)}$ as $(k_1+\frac{1}{2},\dots,k_p+\frac{1}{2}\mid r_1,\dots,r_q)$. Then we have
\[r_j=\mathrm{card}\{i\mid q-k_i\ge j\}.\]
Thus
\begin{equation}\label{Dw2}
\begin{split}
	\langle \rho_n^{(w)},\rho_c\rangle&=\sum_{i=1}^p(p-i)\left(k_i+\frac{1}{2}\right)+\sum_{i=1}^p\sum_{j=1}^{q-k_i}\left(q-j+\frac{1}{2}\right)\\
	&=\sum_{i=1}^p(p-i)\left(k_i+\frac{1}{2}\right)+\sum_{i=1}^p\frac{q^2-k_i^2}{2}\\
	&=\frac{pq^2}{2}+ \frac{1}{2}\sum_{i=1}^{p-1}i(i+1)-\frac{1}{2}\sum_{i=1}^p(k_i-p+i)^2\\
	&=\frac{pq^2}{2}+\frac{p(p^2-1)}{6}-\frac{1}{2}\sum_{i=1}^p(k_i-p+i)^2.
\end{split}
\end{equation}
On the other hand, by (61) of \cite{DY2}, we have
\begin{equation}\label{w1}
	\{1, 2, \dots, q-k_i\}=\{j:e_{p+j}-e_i\in w\Delta^+(\mfg,\mft)\}
\end{equation}
and
\begin{equation}\label{w2}
	\{1, 2, \dots, p-r_j\}=\{i:e_{i}-e_{p+j}\in w\Delta^+(\mfg,\mft)\}.
\end{equation}

Now let $d=p-q$ and
\[C=\sum_{r=1}^{d-1}r^2+\max\{p(2q-1),(p-1)(2q+1)\}.\]
As noted in the introduction, the sum  $\sum_{r=1}^{d-1}r^2$ is understood to be zero when $d\le 0$. 
The "max" term in $C$ can be replaced by 
$$
(p-1)(2q+1)+\max\{1-2d, 0\}
\mbox{ or }
p(2q-1)+\max\{2d-1, 0\}.
$$
Note that
\[-2C=B_0-4\left(\frac{pq^2}{2}+\frac{p(p^2-1)}{6}\right)-\min\left\{\frac{p}{2},\frac{2q+1}{4}\right\}.\]
Then by \eqref{Dw1} and \eqref{Dw2},
\[B_0-D_w(\mu)\ge 2F_w(\mu)+2\sum_{i=1}^p(k_i-p+i)^2-2C.\]
We will show that
\begin{equation}\label{eqn-C}
\F_w(\mu)+\sum_{i=1}^p(k_i-p+i)^2\ge C
\end{equation}
for $a_1>q$, $b_1\ge p-\frac{1}{2}$ and thereby completing the proof of Theorem \ref{thm-main}.

\begin{prop}
	Take $\mu$ as in \eqref{mu}.  Choose $w\in \Delta(\mfg,\mft)^1$ such that $\mu+2\rho_c\in \C_w$. Write $\rho_n^{(w)}$ as $(k_1+\frac{1}{2},\dots,k_p+\frac{1}{2}\mid r_1,\dots,r_q)$. If $a_1>q$ and $b_1\ge p-\frac{1}{2}$, then \eqref{eqn-C} holds.
\end{prop}

\begin{proof}
	 Note that $\langle \mu,\alpha\rangle$ and $\langle\rho_n^{(w)},\alpha\rangle$ are nonnegative for all $\alpha\in \Delta^+(\mfk,\mft)$. Then after discarding the roots  $e_i-e_j \,(i<j)$, $e_i+e_j \,(1<i<j\le p)$ in \eqref{Fw}, we obtain
	\begin{equation}\label{Fw1}
		F_w(\mu)\ge \sum_{i=2}^p\min\{a_1+a_i, k_1+k_i+1\}+\sum_{1\le s<t\le q}\min\{b_s+b_t, r_s+r_t\}+\sum_{j=1}^q\min\{b_j, r_j\}.
	\end{equation}
		Since $k_i\le q$, then $p-i-k_i\ge d-i$ for $1\le i\le d-1$. We have
	\begin{equation}\label{C-1}
		\sum_{i=1}^p (k_i-p-i)^2\ge\sum_{i=1}^{d-1}(d-i)^2=\sum_{r=1}^{d-1}r^2.
	\end{equation}
	
	Firstly, assume that $k_1=q$. Let $l$ denote the number of indices $i$ such that $k_i=q$. Then $r_1=p-l$ and \eqref{w2} shows that
	\[\langle\mu+2\rho_c, e_i-e_{p+1}\rangle\ge 0\]
	for all $1\le i\le l$. The above inequality is equivalent to
	\[a_i\ge b_1+2q-1-2(p-i)\ge q+\frac{1}{2}.\]
	This shows that
	\[\min\{a_1+a_i, k_1+k_i+1\}=2q+1\]
	for $2\le i\le l$ and we always have
	\[\min\{a_1+a_i, k_1+k_i+1\}\ge q+1\]
	for $i>l$. On the other hand, since $r_1=p-l<b_1$, consider only the summands whose indices contain $1$, we obtain that
	\[\sum_{1\le s<t\le q}\min\{b_s+b_t, r_s+r_t\}+\sum_{j=1}^q\min\{b_j, r_j\}\ge (q-1)(p-l)+(p-l)=q(p-l).\]
	Thus \eqref{Fw1} shows that
	\[F_w(\mu)\ge (2q+1)(l-1)+(q+1)(p-l)+q(p-l)=(p-1)(2q+1).\]
	On the other hand, we have
	\[\sum_{i=1}^p (k_i-p-i)^2\ge\sum_{r=1}^{d-1}r^2+\max\{1-2d,0\}.\]
	Indeed, if $\max\{1-2d,0\}> 0$, then $d\le 0$ and thus $\sum_{r=1}^{d-1}r^2=0$. We have
	\[\sum_{i=1}^p (k_i-p-i)^2\ge(k_1-p-1)^2=(d-1)^2\ge 1-2d.\]
	If $\max\{1-2d,0\}=0$, this is exactly \eqref{C-1}. Thus
	\[\F_w(\mu)+\sum_{i=1}^p(k_i-p+i)^2\ge (p-1)(2q+1)+\sum_{r=1}^{d-1}r^2+\max\{1-2d,0\}=C.\]
	This finishes the proof for $k_1=q$.
	
	Now assume that $k_1<q$. This implies that $r_1=p$. Since $k_i\le q-1$, then
	\begin{equation}\label{C-2}
	\begin{split}
		\sum_{i=1}^p(k_i-p+i)^2&=(k_1-p+1)^2+k_p^2+\sum_{i=2}^{p-2}(k_i-p+i)^2\\
		&\ge (k_1-p+1)^2+k_p^2+\sum_{r=1}^{d-1}r^2,
	\end{split}
	\end{equation}
	where the late step uses $p-i-k_i\ge d-i+1$ for $2\le i\le d+1$.
	By \eqref{w1},
	\[\langle\mu+2\rho_c, e_{p+j}-e_1\rangle\ge 0\]
	for $1\le j\le q-k_1$. The above equality is equivalent to
	\begin{equation}\label{w3}
		b_j+2(q-j)+1\ge a_1+2(p-1)\ge q+2p-\frac{3}{2}.
	\end{equation}
	Note that
	\[\min\{a_1+a_i,k_1+k_i+1\}\ge k_1+1,\quad \min\{b_1+b_t,r_1+r_t\}\ge p.\]
	Then consider only the summands whose indices contain $1$, we have
	\[\sum_{i=2}^p\min\{a_1+a_i, k_1+k_i+1\}+\sum_{1\le s<t\le q}\min\{b_s+b_t, r_s+r_t\}\ge (p-1)(k_1+1)+p(q-1).\]
	On the other hand, it is easy to see that
	\[\min\{b_1, r_1\}+\min\{b_q, r_q\}+k_p^2\ge p.\]
	Thus by \eqref{Fw1} and \eqref{C-2},
	\begin{align*}
		F_w(\mu)+\sum_{i=1}^p(k_i-p+i)^2&\ge F_w(\mu)+(k_1-p+1)^2+k_p^2+\sum_{r=1}^{d-1}r^2\\
		&\ge (p-1)(k_1+1)+pq+(k_1-p+1)^2+\sum_{r=1}^{d-1}r^2.
	\end{align*}
	This implies that
	\begin{equation}\label{w4}
		F_w(\mu)+\sum_{i=1}^p(k_i-p+i)^2-\sum_{r=1}^{d-1}r^2-p(2q-1)\ge(k_1-p+1)(k_1-p)-p(q-k_1-1).
	\end{equation}
	We will prove \eqref{eqn-C} according to the following two cases.
	
	(i) $b_{q-k_1}=0$. Take $j=q-k_1$ in \eqref{w3}, we obtain
	\begin{equation}\label{q-k1}
		q\ge 2p+2(q-k_1)-2.
	\end{equation}
	Thus $k_1-p\ge p+(q-k_1)-2$. Therefore,
	\[(k_1-p+1)(k_1-p)\ge p(q-k_1-1).\]
	By \eqref{w4},
	\[F_w(\mu)+\sum_{i=1}^p(k_i-p+i)^2\ge p(2q-1)+\sum_{r=1}^{d-1}r^2.\]
	Note that \eqref{q-k1} implies $d<0$, which in turn shows that $C=p(2q-1)+\sum_{r=1}^{d-1}r^2$. Thus
		\[F_w(\mu)+\sum_{i=1}^p(k_i-p+i)^2\ge C.\]
	
	(ii) $b_{q-k_1}>0$. We \emph{claim} that the omitted summands in $F_w(\mu)$ have sum at least $p(q-k_1-1)$. From this and \eqref{w4}, we have
	\begin{equation}\label{C-3}
		F_w(\mu)+\sum_{i=1}^p(k_i-p+i)^2\ge p(2q-1)+\sum_{r=1}^{d-1}r^2+(k_1-p+1)(k_1-p).
	\end{equation}
	
	Indeed, let $h$ be the smallest index $j$ satisfying $b_j\ge p$ and $2\le j\le q-k_1$. If $h=q-k_1$, then
	\[\sum_{j=2}^{q-k_1}\min\{b_j, r_j\}=\sum_{j=2}^{q-k_1}\min\{b_j,p\}=p(q-k_1-1)\]
	and we are done. Now assume that $h<q-k_1$. In this case, we have $b_h\le p-\frac{1}{2}$. By \eqref{w3},
	\[q-h\ge p+h-2\ge p.\]
	Note that $b_h\ge b_{q-k_1}\ge 1$, then
	\[\sum_{h<t\le q}\min\{b_h+b_t, p+r_t\}\ge (q-h)\ge p.\]
	This proves the claim. Now $k_1-p=-d-(q-k_1)\le -d-1$ shows that
	\[(k_1-p-1)(k_1-p)\ge d(d+1)\ge 2d-1\]
	for $d\ge 1$. Thus $(k_1-p-1)(k_1-p)\ge 0$ implies that
	\[(k_1-p-1)(k_1-p)\ge \max\{2d-1,0\}.\]
	Therefore, by \eqref{C-3}, we have
	\[F_w(\mu)+\sum_{i=1}^p(k_i-p+i)^2\ge p(2q-1)+\sum_{r=1}^{d-1}r^2+\max\{2d-1\}=C.\]
     This finishes the proof for $k_1<q$.
\end{proof}
\medskip

\centerline{\scshape Funding}
Dong is supported by the National Natural Science Foundation of China (no.
12671039).

\end{document}